\documentclass[11pt]{amsart}

\usepackage[backref,bookmarks=false]{hyperref}

\usepackage{color}

\newcommand{\red }{\color{red}}

\numberwithin{equation}{section}

\newcommand{\ben}{\begin{enumerate}}
\newcommand{\een}{\end{enumerate}}

\newcommand{\bea}{\begin{eqnarray}}
\newcommand{\ba}{\begin{array}}
\newcommand{\bean}{\begin{eqnarray*}}
\newcommand{\ea}{\end{array}}
\newcommand{\eea}{\end{eqnarray}}
\newcommand{\eean}{\end{eqnarray*}}
\newcommand{\beq}{\begin{equation}}
\newcommand{\eeq}{\end{equation}}
\newcommand{\bthm}{\begin{thm}}
\newcommand{\ethm}{\end{thm}}
\newcommand{\blem}{\begin{lem}}
\newcommand{\elem}{\end{lem}}
\newcommand{\bprop}{\begin{prop}}
\newcommand{\eprop}{\end{prop}}
\newcommand{\bcor}{\begin{cor}}
\newcommand{\ecor}{\end{cor}}
\newcommand{\bdfn}{\begin{dfn}}
\newcommand{\edfn}{\end{dfn}}
\newcommand{\brem}{\begin{rem}}
\newcommand{\erem}{\end{rem}}
\newcommand{\bpf}{\begin{proof}}
\newcommand{\epf}{\end{proof}}
\newcommand{\bfact}{\begin{fact}}
\newcommand{\efact}{\end{fact}}
\newcommand{\bobs}{\begin{obs}}
\newcommand{\eobs}{\end{obs}}
\newcommand{\bexam}{\begin{exam}}
\newcommand{\eexam}{\end{exam}}
\newcommand{\bclaim}{\begin{claim}}
\newcommand{\eclaim}{\end{claim}}

\newtheorem{thm}{Theorem}[section]
\newtheorem{prop}[thm]{Proposition}
\newtheorem{lem}[thm]{Lemma}

\newtheorem{cor}[thm]{Corollary}
\newtheorem{dfn}[thm]{Definition}
\newtheorem{rem}[thm]{Remark}
\newtheorem{fact}[thm]{Fact}
\newtheorem{claim}[thm]{Claim}
\newtheorem{obs}[thm]{Observation}
\newtheorem{exam}[thm]{Example}

\newtheorem{condition}{Condition}
\newtheorem*{condition'}{Condition 2'}

 \newtheoremstyle{claimstyle}%
   {}
   {}
   {\normalfont}
   {}
   {\itshape}
   {.}
   { }
   {\thmnote{#3}}

\theoremstyle{claimstyle}

\alph{enumii} \roman{enumiii}

\def\cA{\mathcal A}                    
             \def\cF{\mathcal F}       
\def\cL{{\mathcal L}}                   
                    \def\cJ{\mathcal J}   \newcommand{\J}{\mathcal{J}}
                           \def\cE{\mathcal E}

                \def\Z{{\mathbb Z}}      \def\R{{\mathbb R}}
\def\C{{\mathbb C}}

\newcommand{\cbar}{\hat{{\mathbb C}} }

\def\a{\alpha}                             \def\d{\delta}
\def\De{\Delta}               \def\e{\varepsilon}           
\def\g{\gamma}                \def\Ga{\Gamma}           \def\l{\lambda} \def\la{\lambda}
\def\La{\Lambda}              \def\om{\omega}           
               \def\sg{\sigma}
               \def\th{\theta}           
\def\ka{\kappa}

\newcommand{\ep}{\varepsilon}
\newcommand{\ph}{\varphi}
\newcommand{\al}{\alpha}
\newcommand{\ga}{\gamma}

\def\1{1\!\!1}

\def\and{\text{ and }}

     \def\HD{\text{{\rm HD}}}

              \def\bu{\bigcup}
\def\({\bigl(}                \def\){\bigr)}
\def\lt{\left}                \def\rt{\right}

\def\ld{\ldots}                        \def\^{\tilde}

\def\sbt{\subset}             \def\spt{\supset}

           \def\downto{\searrow}
\def\sp{\medskip}             \def\fr{\noindent}        
\def\ov{\overline}            
\def\ess{{\rm ess}}
\def\fr{\noindent}

\def\om{\omega}

\def\arg{\text{arg}}

\def\EP{\mathcal E{\text{{\rm P}}}}
\def\D{{\mathbb D}}

\def\${$ \displaystyle }

\def\shift{\theta}

\newcommand{\pf}{{\mathcal{L}}}

\newcommand{\npf}{\mathcal{\hat L}}

\newcommand{\jul}{\mathcal J}

\newcommand{\sT}{\mathring T}

\def\rad{\J_r}

\def\mirror{\Upsilon}
\def\analytic{\cA} \def\ranalytic{\cA_\R}
\def\nlip{Lip_n}

\def\whatever{admissible }

\def\norm{{\bf |}}

\begin{document}

\title[Real Analyticity for random dynamics]{Real Analyticity for random dynamics of transcendental functions}


\author{Volker Mayer}
\address{Volker Mayer, Universit\'e de Lille I, UFR de
  Math\'ematiques, UMR 8524 du CNRS, 59655 Villeneuve d'Ascq Cedex,
  France} \email{volker.mayer@math.univ-lille1.fr \newline
  \hspace*{0.42cm} \it Web: \rm math.univ-lille1.fr/$\sim$mayer}

\author{Mariusz Urba\'nski}
\address{Mariusz Urba\'nski, Department of
  Mathematics, University of North Texas, Denton, TX 76203-1430, USA}
\email{urbanski@unt.edu \newline \hspace*{0.42cm} \it Web: \rm
  www.math.unt.edu/$\sim$urbanski}
  
  \author[\sc Anna ZDUNIK]{\sc Anna Zdunik}
\address{Anna Zdunik, Institute of Mathematics, University of Warsaw,
ul. Banacha 2, 02-097 Warszawa, Poland}
\email{A.Zdunik\@@mimuw.edu.pl}

\date{\today} \subjclass{111}

\thanks{A. Zdunik was supported in part by  NCN grant 2014/13/B/ST1/04551.
The research of M. Urba\'nski was supported in part by the NSF Grant DMS 1361677.
This work was also supported in part by the Labex CEMPI  (ANR-11-LABX-0007-01).}


\begin{abstract}  Analyticity results of expected pressure and invariant densities in the context of random dynamics of transcendental functions are established. These are obtained by a refinement of work by Rugh \cite{Rug08} leading to a simple approach to analyticity. We work under very mild dynamical assumptions. Just the iterates of the Perron-Frobenius operator are assumed to converge. 

We also provide  Bowen's formula expressing the almost sure Hausdorff dimension of the radial fiberwise Julia sets in terms of the zero of an expected pressure function. 
Our main application establishes real analyticity for the variation of this dimension for suitable hyperbolic random systems of entire or meromorphic functions.
 \end{abstract}

\maketitle

\section{Introduction}\label{intro} Answering a conjecture of Sullivan, Ruelle \cite{Ru82} showed  for hyperbolic rational functions  that the Hausdorff dimension of the Julia sets does depend analytically on the map and gave a local formula for perturbations of the map $z\mapsto z^2$. Since then, there where several results of this type in various contexts and also different methods of proof. The monograph \cite{Zins2000} treats the local formula and analyticity
 has been obtained, for example, in \cite{VW96} for complex Henon mappings of $\C^2$, in \cite{Po15} for basic sets of surface diffeomorphisms. In the context of entire and meromorphic functions, the first result was obtained in \cite{UZ04},
 further development appeared in  \cite{MyUrb08, MUmemoirs} and \cite{SU14}.
 
 Whereas the latter papers use holomorphic motions, Rugh \cite{Rug08} introduces the method of positive cones and complex cones which allowed him to 
 extend analyticity results to random dynamics of repellers. 
 The present paper refines  Rugh's approach, avoids complex cones, and allows us to get analyticity results for random dynamics of transcendental entire and meromorphic functions. The following is a particular case of our general result Theorem \ref{thm main}.

\bthm\label{thm intro 1}
Let $f_\eta (z)=\eta e^z$  and let $a\in (\frac{1}{3e}, \frac{2}{3e})$
and $0<r< r_{max}$, $r_{max}>0$. Suppose that $\eta_1,\eta_2,..$ are i.i.d. random variables uniformly distributed in $\D (a,r)$. Let $J_{\eta_1,\eta_2,...}$ denote the Julia set of the sequence of compositions $\(f_{\eta_n}\circ f_{\eta_{n-1}}\circ\ldots\circ f_{\eta_2}\circ f_{\eta_1}:\C\to\C\)_{n=1}^\infty$
and let 
$$
\rad (\eta_1,\eta_2,...)= \{z\in J_{\eta_1,\eta_2,...}  : \; \liminf_{n\to\infty} |f_{\eta_n}\circ ...\circ f_{\eta_1}(z)|<+\infty\}
$$
be the radial Julia set of $(f_{\eta_n}\circ ...\circ f_{\eta_1})_{n=1}^\infty$. 
Then, the Hausdorff dimension of $\rad (\eta_1,\eta_2,...)$ is almost surely constant and depends real-analytically on the parameters $(a,r)$ provided that $r_{max}$ is sufficiently small.
\ethm

The common point in all papers on this topic is the fact that the Hausdorff dimension of Julia sets can be expressed in terms of the zero of a pressure function. This fact goes back to  \cite{Bow79} and is now called Bowen's Formula. This formula also has been generalized in many contexts and we also provide one (Theorem \ref{thm BF}). We would like to mention that the zero of the involved (expected in the random case) pressure  does not really detect the dimension of the whole Julia set but the dimension of its subset consisting of all radial points. In fact, in the case of hyperbolic rational functions the radial Julia set and the Julia set itself coincide. However, for transcendental functions, especially for entire functions, there is a definite difference between these sets. McMullen \cite{McM87} showed that the Julia set of sine or exponential functions is always maximal equal to two whereas for such hyperbolic functions the dimension of the radial Julia set, which is often called hyperbolic dimension,  is never equal to two \cite{St96, UZ03}.

\medskip

The formulation of Theorem \ref{thm intro 1} has been chosen deliberately in analogy with Example 1.2 in 
\cite{Rug08} since our present work stems from Rugh's papers \cite{Rug08, Rug10}.
However, we were not able to apply directly his machinery. Instead we worked out a refinement
of Rugh's elegant approach to analyticity. In particular, we avoid any use of Hilbert's distance in positive cones and complex cones. Instead we provide a quite simple and direct calculation (see Proposition \ref{34}).
The outcome, besides the results concerning random transcendental dynamics, provides an elementary and general tool.
In short, it says that if the thermodynamic formalism holds and if the normalized iterated transfer operator converges with a uniform speed, then real analyticity holds. Let us explain this now in more detail.

We consider arbitrary analytic families of holomorphic functions $f_{j,\l}$, $j\in \Z$, having the following properties. There exists an open set $U\subset \cbar$ and $\d_0\in (0, 1/4)$ such that, for all $w\in  U$, $j\in \Z$ and $n\geq 1$, every inverse branch $g$ of the non-autonomous composition 
$$f_{j,\l}^n:= f_{j+n-1,\l}\circ ...\circ f_{j+1,\l}\circ f_{j,\l}$$
 exists on $\D(w,2\d_0)$, maps $\D(w,\d_0)$ inside $U$ and satisfies $|g'|\leq \g_n^{-1}$ on this disk. Here $(\ga_n)_n$ is any sequence with $\lim_{n\to\infty}\ga_n=\infty$. As for specific examples, the reader my have in mind rational functions, functions associated to finite or infinite iterated function systems or transcendental functions.
In such a setting the thermodynamical formalism including a Ruelle-Perron-Frobenius Theorem usually holds (see for example \cite{Rue78}, \cite{Zins2000}, \cite{PUbook}, \cite{MauldinUrb03}, \cite{Kif08}, \cite{MSUspringer}, \cite{MUmemoirs} and \cite{MyUrb2014}): 
$$\pf_{j,\l , t} g (w) = \sum_{f_{j,\l}(z)=w}|f_{j,\l}'(z)|^{-t}\left(\frac{1+|z|^2}{1+|w|^2}\right)^{-\tau \frac{t}{2}}\!\!\!g(z)
\quad , \quad g\in C_b^0 (U)\; ,$$
defines a bounded operator on the space of bounded continuous functions $C_b^0 (U)$ equipped with the sup-norm
such that, for every $j\in \Z$,
\ben
\item[-] there exists probability measures $\nu_{j,\l ,t}$ and reals $P_{j}(\l, t)$ such that 
\beq\label{conf measures}\pf_{j,\l , t}^*\nu_{j+1,\l ,t} = e^{P_{j}(\l, t)}\nu_{j,\l ,t}\eeq
\item[-] and that there exist functions $\hat \rho _{j,\l ,t}\in C_b^0 (U)$ such that $\npf _{j,\l ,t}\hat \rho _{j,\l ,t}=\hat \rho _{j+1,\l ,t}$ where  $\npf _{j,\l ,t}=e^{-P_{j}(\l, t)} \pf _{j,\l ,t}$ is the normalized operator.
\een
The functions $\hat \rho _{j,\l ,t}$, called invariat densities give rise to an invariant family of measures $\mu_{j,\l ,t}$, defined as  $d\mu_{j,t,\l}=  \hat \rho _{j,\l ,t}d\nu_{j,\l,t}$. This family is invariant in a sense that
$(f_{j,\l})_*(\mu_{j,\l,t})=\mu_{j+1,\l,t}$.

 In here, $t$ belongs to an interval $I$ of positive reals and  $\tau \geq 0$. When $\tau =0$ then the above operators are just the usual geometric transfer operators used, for example, for polynomials or iterated function systems.
 For the infinite to one transcendental functions we have to use  the additional coboundary factor with some well chosen $\tau>0$.
 
 In such a setting the iterated normalized operators are uniformly bounded, i.e. there exists $M<\infty$ such that
 \beq\label{42}
 \| \npf _{j,\l ,t}^n \|_\infty  \leq M \quad \text{for all} \quad j\in \Z\, ,\l\in \La \text{ and }  t\in I
 \eeq
 where $\npf _{j,\l ,t}^n = \npf _{j+n-1,\l ,t} \circ ...\circ \npf _{j,\l ,t}$ (see  \cite{MyUrb08} and \cite{MyUrb2014} for the case of transcendental functions). Also, the densities satisfy the following positivity condition as soon as the dynamical system is mixing (see for example Lemma 5.5 in \cite{MyUrb2014}): there exists $z_0\in U$ and $a>0$ such that
\beq \label{43}
\hat \rho _{j,\l , t} (z_0) \geq a  \quad \text{for all} \quad j\in \Z\, ,\l\in \La \text{ and }  t\in I\,.
\eeq
We use a bounded deformation property. It is formulated in Definition \ref{15} and gives a uniform control of the variation 
of local inverse branches. 
Finally, to $\eta >0$ we associate the space of Lipschitz functions $Lip(U,\eta )$
which is the space of bounded functions $g:U\to \R$ such that 
\beq\label{Lip eta}
Lip(g, \eta)=\sup \left \{  \frac{|g(z_1)-g(z_2)|}{|z_1-z_2|} \; ; \; z_1,z_2\in U \; , \; 0<|z_1-z_2|<\eta\right\}<\infty \, .
\eeq
This space is equipped with the norm $\| g\|_{Lip , \eta} = \|g\|_\infty + Lip(g , \eta)$.

\smallskip

\bthm \label{thm intro 2}
Suppose that $f_{j,\l}$ are of bounded deformation and that  the above thermodynamical formalism holds, in particular with \eqref{42} and \eqref{43}. 
 Suppose that the iterated normalized operators have uniform speed:  for every $\eta >0$ there exists $\om_n\to 0$ such that 
\beq\label{uniform speed}
\| \hat \pf^n_{j, \l,t} g -\nu_{j,\l,t} (g) \hat \rho_{j+n,\l,t} \|_{\infty} \leq \om_n \|g\|_{Lip , \eta} \quad \text{for every} \;\; \;  n\geq 1
\; , \;\; j\in \Z \; \;\text{and}
\eeq
every $ g\in Lip (U , \eta)$. Let $z_0\in U$.
Then, for every $j\in \Z$, the function
$$
(\l,t,z) \mapsto \rho_{j,\l ,t}(z)= \frac{\hat  \rho_{j,\l ,t}(z)}{\hat  \rho_{j,\l ,t} (z_0)} 
$$ is real analytic.
\ethm

Theorem \ref{thm intro 2} will be a consequence of Theorem \ref{64}, Theorem \ref{thm analyticity} is its  random analogue. 
All these results concern real analyticity of invariant densities. In fact Theorem \ref{64} proves a stronger version of real analyticity than the one in Theorem\ref{thm intro 2}; namely that the mapping
$$
(\l,t)\mapsto \rho_{j,\l ,t}
$$
is real-analytic, where $\rho_{j,\l ,t}$ is considered as a member of an appropriate natural Banach space.

As it is explained in Remark \ref{68}, Theorem \ref{thm analyticity}  could also include real analyticity of expected pressure. We worked this out in detail in the case of  random transcendental dynamics and the cumulating result including real analyticity of the hyperbolic dimension is Theorem \ref{thm main}.

\section{General setting}\label{General Setting}
We already outlined the setting in the Introduction and present now details.
They will be formulated for the non-autonomous setting since all the sections to follow including 
Section 7 are devoted to non-autonomous dynamics. Random dynamics are the object of Section 8 and 9.
We denote by $D_z=\D(z,\d )$ the Euclidean disk of radius $\d$ centered at $z\in \C$. Suppose given
$$\text{an open set 
$\; U\subset \C\; $, $\; 0<\d <\d_0<\frac14\; $
and a sequence $\;\ga_n\to \infty$.}$$
 
 \noindent
For $j\in \Z$, we suppose that $ f_j$ is a holomorphic function defined on some open set $ V_{f_j}\subset \C$ with range in
 $\cbar$  such that the following holds for  every $j\in \Z$ and $n\geq 1$:
 the composition $f_j^n= f_{j+n-1}\circ ...\circ f_j$
is defined on some domain such that the range contains the euclidean $2\d_0$--neighborhood of $U$ and such that
for every $w\in U$ 
 every inverse branch $g$ of $f_j^n$ is well defined on $\D(w, 2\d _0)$ and satisfies
\beq\label{20}
g(\D(w,\d_0))\subset U\quad \text{and}\quad |g' |_{|\D(w,\d_0)} \leq \ga _n ^{-1}\,.
\eeq

As often, replacing the functions by some of their iterates, we can assume that $\ga_n>1$ for all $n\geq 1$.

\bexam
The reader may have in mind the following examples:
\ben
\item[-] $f_{j,\l} (z) = z^2+\l c _j$ where $\l \in \D(0,1)$ and $|c_j |< \frac18$ or other suitable perturbations of hyperbolic rational functions.
\item[-] Functions arising from (finite or infinite alphabet) conformal iterated functions systems.
\item[-] Families of transcendental functions such as the exponential family in Theorem \ref{thm intro 1} and all the examples treated in \cite{MUmemoirs, MyUrb2014}.
\een
\eexam

From the above definition follows that every function $f_j$ has the set $U$ in its range $V_{f_j}$ and that
$f_j^{-1} (\ov U)\subset U$. As a motivation for our non-autonomous setting, note that the radial Julia set of a hyperbolic meromorphic function $f:\C\to\hat\C$ (see \cite{MyUrb08} and \cite{MUmemoirs} for a precises definition of this concept) is
$$
\rad (f) =\left\{ z\in \bigcap_{n>0}f^{-n}(U)\; : \;\; \liminf_{n\to\infty} |f^n(z)|<+\infty\right\}\,,
$$
where $U$ here is a sufficiently small neighborhood of the Julia set $J(f)$.
The straightforward adaption of this definition to the non-autonomous case is the following:
\beq\label{51}
\rad (f_j,f_{j+1}, ... ) =\left\{ z\in \bigcap_{n>0}(f_{j}^n)^{-1}(U)\; : \;\; \liminf_{n\to\infty} |f_j^n(z)|<+\infty\right\}\,,
\eeq
where $U$ is now as above, for ex. in \eqref{20}.
Notice that these radial Julia sets coincide with the usual Julia sets (of the same sequence) as soon as the open set $U$ is bounded. This is the case for rational functions (after an appropriate change of coordinates) and for iterated function systems. Unbounded sets $U$ and radial Julia sets are necessary for transcendental dynamics.

Our results concern holomorphic families of functions. Let $\La$ be the corresponding parameter space.
Without loss of generality, we may assume that $\La$ is one-dimensional and, the results being of local nature, we can restrict to the case where $\La= \D (\l_0 , r)$ is an open disk in $\C$ having arbitrarily
small radius $r>0$.

\bdfn \label{21} $\cF_{\La} = \{f_{j,\l}\; , \; j\in \Z\; and \; \l\in \La\}$ is called a non-autonomous holomorphic family if, for every $j\in \Z$, $f_{j,\l }$ depends holomorphically on $\l\in \La$. This precisely means the following for every $j\in \Z$: if $V_{f_{j,\l}}$ is the domain of $f_{j,\l}$ then
$\Ga_j:=\bigcup_{\l\in \La}\{\l \}\times V_{f_{j,\l}}$ is an open subset  of $\C^2$ and the map $ (\l,z)\mapsto f_{j,\l}(z)$ is holomorphic on $\Ga_j$.
\edfn

\section{Pairings and bounded deformation}
Let $\cF_{\La}=\{f_{j,\l} \; ; \; j\in \Z \; , \; \l\in \La\}$ be a non-autonomous holomorphic family.
We are interested in the solutions of the equation $f_{j,\l} ^n (z)=w$.
A direct application of the implicit function theorem along with analytic continuation and \eqref{20}
gives the following observation.

\bfact\label{12}
If $\l ' \in \La$, $w'\in U$ and $z'\in f_{j,\l '}^{-n} (w')$ are given, then there exists a unique holomorphic function
$$
\begin{array}{clc}
\La \times  \D(w' , \d _0) &\to &U\\
(\l , w) & \mapsto &z(\l , w)
\end{array}
$$
such that $z(\l ', w') =z' $ and $f_{j,\l} ^n (z(\l, w)) =w $ for every $\l \in \La $ and $w\in \D (w', \d _0)$.
\efact

\fr This is simply the proper way of defining an inverse branch $f_{j,\l , *}^{-n}$ of $f_{j,\l}^n$. We will use the inverse branch notation rather than the function $z(\l,w)$. This precisely means that $f_{j,\l , *}^{-n}$ is a choice of inverse branch defined by 
a function $z$ given by Fact \ref{12}:
$$
f_{j,\l , *}^{-n} (w) = z(\l , w) \quad , \quad \l \in \La \; , \; w\in \D (w', \d _0)\,.
$$

We can now introduce the notion of pairings used in the sequel. Let us recall that $0<\d\leq \d_0$. The number $\d$ will be specified later on
in \eqref{small delta}.

\bdfn\label{13} $(w_1,w_2)$ is a $0$--pairing if $w_1\in U $ or $w_2\in U$ and if $|w_1-w_2|<\d$.
For $n\geq 1$, $(z_1,z_2)$ is called $n$--pairing if there exists a $0$--pairing $(w_1,w_2)$, $j\in \Z$, parameters $\l_1,\l_2\in \La$ and a choice of inverse branch $f_{j,\l , *}^{-n}$ such that 
$$ z_1= f_{j,\l_1 , *}^{-n}(w_1) \quad \text{and}\quad  z_2= f_{j,\l_2 , *}^{-n}(w_2) \; .$$
\edfn

\smallskip

The following concept of bounded deformation has already been used in \cite{MUmemoirs} but without the condition \eqref{bdd def 2}.
This was so since for dynamically regular transcendental functions this second condition automatically is satisfied
(see Lemma \ref{lemma 9.7}). It is also possible to relax this second condition in the setting of conformal infinite iterated function systems as it has been done \cite{SU14}.

\bdfn\label{15}
The family $\cF_{\La}$ is of bounded deformation if there exists $A, D<\infty$ such that for every $j\in \Z$ and for every choice of inverse branch
$f_{j,\l, *}^{-1}$  we have 
\beq\label{bdd def 1}
\left|\frac{\partial f_{j,\l,*}^{-1}}{\partial \l}\right| \leq D \; , \quad  \l\in \La
\quad  \text{ and } 
\eeq
\beq\label{bdd def 2}
\left|  \frac{f_{j,\l_1}'(z_{1})}{f_{j,\l_2}'(z_{2})} \right| =
\left|  \frac{f_{j,\l_1} '(f^{-1}_{j,\l_1 , *}(w_1))}{f_{j,\l_2} '(f^{-1}_{j,\l_2 , *}(w_2))} \right|\leq A
\; , \quad \l_1,\l_2\in \La \; , \;  w_1,w_2\in \D (w,\d_0)
\;  .
\eeq
\edfn

Bounded deformation holds for many transcendental families
and especially for $f_\l (z) = \l e^z$ (see \cite{MUmemoirs}).
Notice that \eqref{bdd def 1} is equivalent to the fact that 
$
\left|\frac{\partial f_{j,\l}}{\partial \l}\right| \leq D \left| f_{j,\l}' \right|
$.
This condition is  automatically satisfied for all rational functions and for functions associated to finite iterated function systems subject to possible shrinking of the parameter space. Also, for all systems with compact phase space such as infinite iterated function systems one can use the theory of holomorphic motions in order to show that \eqref{bdd def 1} holds for free. So, the bounded deformation condition is mainly instrumental in the case of transcendental, and especially entire, functions.

\sp Remember that the expanding constant $\g_1 >1$. This allows us to fix a constant 
\beq \label{kappa}
\kappa \in (\g_1^{-1} , 1)\, .
\eeq

\blem\label{16}
If  $(f_{j,\l})$ satisfies \eqref{bdd def 1}  
then there exists a (sufficiently small) choice of $diam (\La )$ (depending on $\d$) such that every $1$--pairing  $(z_{1}, z_{2})
$
satisfies $|z_{1}-z_{2}|<\ka\d$. 
\elem

\brem \label{rem new 1}
Lemma \ref{16} implies that every $1$--pairing is a $0$--pairing and, inductively, 
that every $n$--pairing is a $k$--pairing for all $0\leq k <n$.
\erem

\bpf
Let a $1$--pairing be given by $z_{i}=f_{j,\l_i,*}^{-1}(w_i)$, $i=1,2$, and denote $z'_{2}=f_{j,\l_2,*}^{-1}(w_1)$. 
The condition \eqref{bdd def 1} implies that
$$ \left|z_{1}-z'_{2}  \right|= \left|f_{j,\l_1,*}^{-1}(w_1)-f_{j,\l_2,*}^{-1}(w_1)  \right|\leq D \; diam( \La )\,.
$$
On the other hand, $|z'_{2} -z_{2} |<\ga_1^{-1}\d$.
Therefore, $|z_{1} -z_{2} |< \d \ga_1^{-1} + D diam(\La )$ and
it suffices to take $diam(\La ) < \d (\ka-\ga_1^{-1})/D$.
\epf

\
 A further consequence of bounded deformation, this time of condition \eqref{bdd def 2}, is the following.

\blem\label{44}
 There exists  a constant  $\tilde A<\infty$ independent of $\d\in (0,\d_0)$
 such that, for every $j\in \Z$ and every $1$--pairing 
 $ (z_1,z_2)= (f_{j,\l_1,*}^{-1}(w_1), f_{j,\l_2,*}^{-1}(w_2))$,
\beq \label{41}
\left| \arg \left(  \frac{f_{j,\l_1} '(z_{1})}{f_{j,\l_2} '(z_{2})} \right)\right|=
\left| \arg \left(  \frac{f_{j,\l_1} '(f^{-1}_{j,\l_1 , *}(w_1))}{f_{j,\l_2} '(f^{-1}_{j,\l_2 , *}(w_2))} \right)\right|\leq \tilde A\,
\eeq
provided the parameters $\l_1,\l_2\in \D(\l_0, r/2)$.
 In here, the argument is well defined and understood to be the principal choice, i.e. $\arg (1)=0$.
\elem

\bpf
By Koebe's distortion theorem (see for ex. Theorem 2.7 in \cite{McMullen94}) it suffices to consider pairings for which 
$f_{j,\l_1}(z_{1})=f_{j,\l_2}(z_{2})=w\,$
or, in terms of inverse branches, that $z_{i}=f_{j,\l_i,*}^{-1} (w)$, $i=1,2$. Consider then the function
$$\ph (\l ) = \frac{f_{j,\l} '(f^{-1}_{j,\l , *}(w))}{f_{j,\l_0} '(f^{-1}_{j,\l_0 , *}(w))}\quad , \quad \l \in \La =\D(\l_0 ,r)\,.
$$
It has the properties $\ph (\l_0)=1$ and $A^{-1}\le |\ph |\leq A$ by \eqref{bdd def 2}.
The set of all holomorphic functions having these properties is compact which implies the estimate 
\eqref{41}.
\epf

\
In the rest of this paper we suppose that $r=diam\La /2$ is chosen such that the conclusion of Lemma \ref{16} holds
as well as
\eqref{41} for every $1$--pairing, i.e. \eqref{41} holds for all parameters $\l_1,\l_2\in \D(\l_0, r)$.

\section{Mirror extension}
One step towards real analyticity is complexification of the transfer operator and its potential.
 There are several possibilities for this but the elegant mirror extension of Rugh  is now most appropriate for us. We use mainly the notation he used in his papers \cite{Rug08,Rug10}.
The mirror of the parameter space $\La$ and the domain $U$ is the set 
\beq \label{11} 
\mirror = \Big\{ (\l_1,  \ov \l_2, w_1, \ov w_2)\; : \;\; \l_1,\l_2\in \La \, , \; (w_1, w_2) \text{ is a $0$--pairing}\Big\}\,.
\eeq
Consider also the $w$--mirror 
$$\mirror _w =\Big\{ (w_1, \ov w_2)\; \; : \;\;(w_1, w_2) \text{ is a $0$--pairing}\Big\}\, .$$
The initial sets $ \La\times U$ and $U$ naturally identify respectively with the diagonals
$$\Delta = \{(\l , \ov \l , w ,\ov w ) \, : \;\; \l \in \La\; , \; w\in U\}\subset \mirror
\quad \text{and} \quad \Delta _w= \{(\om ,\ov \om ) \, : \;\;  \om\in U\}\subset \mirror_w\,.$$

\noindent
Let
$\analytic =  C^\om_b (\mirror _w)$ be the space 
 of functions that are holomorphic  and bounded on  $\mirror_w$. 
 This space will be equipped with the sup-norm defined by 
 $$
 \|h\|_\infty = \sup_{(w_1,\ov w_2)\in \mirror_w} |h(w_1, \ov w_2 )]
 $$
and it makes it a Banach space.
We also need the following notion of Lipschitz variation on $n$--pairings of a function $h:\mirror_w \to \C$:
\beq \label{38}
\nlip (h) = \sup \left\{ \frac{|h(z_1, \ov z_2) - h(z_1, \ov z_1)|}{|z_1-z_2|}\; , \;\; (z_1 , z_2) \text{ $n$--pairing with } z_1\neq z_2
\right\}\,.
\eeq

\blem\label{39}
For every $n\geq 1$ and  $h\in \analytic$ we have  $\nlip (h) \leq \| h\|_\infty /((\kappa-\ga_n^{-1})\d)$, i.e. for every $h\in \analytic$ and every $n$--pairing $(z_1, z_2)$ 
\beq\label{39'}
|h(z_1, \ov z_2) -h(z_1, \ov z_1)|\leq \frac{\| h\|_\infty}{(\ka-\ga_n^{-1})\d} |z_1-z_2|\,.
\eeq
with $\ka$ the constant from \eqref{kappa}.
\elem

\bpf
Let $\sg =\partial \D(z_1, \kappa \d)$. Cauchy's Integral Formula implies
\begin{align*}
|h(z_1, \ov z_2)& -h(z_1, \ov z_1)|\leq \\
&\leq\frac{1}{(2\pi)^2}\int _\sg \int _\sg \left|\frac{h(\xi_1, \ov \xi_2)}{(\xi_1-z_1)(\ov \xi _2-\ov z_2)}-\frac{h(\xi_1, \ov \xi_2)}{(\xi_1-z_1)(\ov \xi _2-\ov z_1)}\right| |d\xi_1 | |d\ov \xi _2|.
\end{align*}
Elementary estimations give  $|\xi_i -z_1|=\kappa\d$ and $|\xi_i -z_2|\geq \d(\ka-\ga_n^{-1})$, $i=1,2$.
The required estimation follows now easily.
\epf

\fr
The space $\analytic$ contains  the relevant subspace
$$ \ranalytic =\left\{ h\in \analytic \; : \;\; h_{| \Delta _w} \in \R\right\}\,.$$
Functions from $ \ranalytic$ are real on the diagonal and can therefore be identified with a subclass of real functions defined on $U$. 
Up to identification, they belong to the space of Lipschitz functions $Lip(U, \eta )$ (see  Introduction) provided $\eta < \ka \d$.

\blem\label{LipschitzLemma}
If $h\in \ranalytic$ then $z\mapsto g(z):=h(z,\ov z)$ belongs to $Lip (U, \ka \d)$ and 
$$\|g\|_{Lip, \ka \d}\leq C \|h\|_\infty$$
where $C= 1+2/( (\sqrt{\ka}-\ka)\d)$.
\elem

\bpf Let $h\in \ranalytic$ and let $z_1,z_2\in U$ with $0<|z_1-z_2|< \ka \d $. Consider $\sg = \partial \D(z_1 , \sqrt{\ka } \d)$ and use exactly the same argument then in the proof of Lemma \ref{39} based on Cauchy's Integral Formula in order to obtain the estimates
$$
|h(z_1, \ov z_2) - h(z_1, \ov z_1) | \leq \frac{1}{(\sqrt{\kappa}-\ka)\d}\|h \|_\infty |z_1-z_2|\,.
$$
The same argument also gives the following symmetric version of this estimaties:
$$
|h(z_1, \ov z_2) -h(z_2, \ov z_2)|\leq  \frac{1}{(\sqrt{\kappa}-\ka)\d}\|h \|_\infty |z_1-z_2|\,.
$$
It suffices now to combine these two estimations in order to complete this proof.
\epf

\subsection{Potentials and extended operator} The potentials under consideration must have two properties:
they must admit holomorphic mirror extensions and have good distortion properties.  We do not treat the most general setting but focus in the following on the most important class of potentials and will see that they have the required properties. So, suppose that $\tau \geq 0$ is fixed, that $I$ is an open interval compactly contained  in $ (0,\infty )$, consider
\beq\label{47}
 \ph_{j,\l , t} (z)= -t\log |f'_{j,\l} (z)| -t\frac{\tau}{2}\log \left(\frac{1+|z|^2}{1+|f_{j,\l} (z)|^2}\right)\quad  \eeq
 and observe that $| f_{j,\l} '|_\tau ^{-t}= e^{\ph_{j,\l , t} }$, $ \l \in \La$ and $ t\in I$,
 where $| f_{j,\l} '|_\tau$ denotes the derivative with respect to the Riemannian conformal metric $|dz|/(1+|z|^2)^{\frac{\tau}{2}}$.
The transfer operator $\pf_j=\pf_{j,\l ,t}$ of the function $f_{j,\l}$ and the potential $\ph_{j,\l , t} $ is defined by
\beq\label{1}
\pf_j g (w) =\sum_{f_{j,\l} (z)=w} e^{\ph_{j,\l ,t} (z)} g(z) =
\sum_{f_{j,\l} (z)=w}  | f_{j,\l} '(z)|_\tau ^{-t}g(z)
\quad , \quad w\in U\, , 
\eeq
where $g\in C_b^0 (U)$ is a continuous bounded function on $U$.
The classical case, particularly when one deals with polynomials or iterated function systems, is when $\tau=0$. For transcendental functions $\tau>0$, i.e. the additional coboundary term $\log (1+|z|^2)-\log (1+|f_{j,\l} (z)|^2)$, is needed since otherwise the transfer operator is not well defined the series defining it being divergent.

The n-th composition of these operators is
 \beq\label{non auto comp pf}
 \cL_j^n = \cL_{j+n-1}\circ...\circ \cL_j\,.
 \eeq
A standard calculation shows that $\pf_j^n$ is the transfer operator as defined in \eqref{1} of the potential
  $$ S_n\ph_j = \sum_{k=0}^{n-1} \ph _{j+k}\circ f_j^k= \sum_{k=0}^{n-1} \ph _{j+k, \l ,t}\circ f_{j,\l}^k\,.$$


\

 The potentials defined in \eqref{47}, often called geometric, admit mirror extensions as we explain now.
 In the following, $\mathcal I$ is a complex neighborhood of $I\subset \R$.
 For $w\in U$, define $Z_w= \La\times \ov \La\times  D_w \times \ov {D_w}$ and notice that  $\mirror \subset \bigcup_{w\in U} Z_w$. From Fact \ref{12} applied with $n=1$ follows that, to every choice of $\l'\in \La$ and $z'\in f_{j,\l'}^{-1}(w)$, there corresponds a choice of inverse branches  $f_{j,\l,*}^{-1}$ defined on $\La\times D_w$. Consider then on $Z_w$ the map
 \beq\label{xw1}
  (\l_1,\ov\l_2, w_1, \ov w_2 )\mapsto (\l_1,\ov\l_2, f_{j,\l_1,*}^{-1} (w_1) , \ov {f_{j,\l_2,*}^{-1}(w_2)} )
  \eeq
and denote its range by $Z^{-1}_{j,w,*}$. Notice that Lemma \ref{16} and \eqref{20} imply 
$$ Z^{-1}_{j,w,*} \subset Z_{w'} \cap (\La\times \ov \La \times U\times \ov U)\quad  \text{ for some $ w'\in U $.}$$
Given the definition of the transfer operator in \eqref{1}, it suffices to extend the potentials to 
\beq \label{17}
\mirror^{-1}\times \mathcal I :=\bigcup_{w , *} Z^{-1}_{j,w,*}\times \mathcal I\subset \mirror \times \mathcal I\,.
\eeq
The extension of  $\ph_{j,\l,t}$ to one of the sets  $Z^{-1}_{j,w,*}\times \mathcal I$ is straightforward. Indeed, let
\beq\label{40}
 \Phi_{j,\l_1,\ov \l_2 , t} (z_1, \ov z_2) =-\frac{t}{2}\log \big( f_{j,\l_1}'(z_1) \ov{f_{j,\l_2}'(z_2) }\big)
 -t\frac{\tau}{2}\log \left(\frac{1+z_1\ov z_2}{1+f_{j,\l_1} (z_1) \ov {f_{j,\l_2} (z_2) }}\right)
\eeq
where $(\l_1 ,\ov\l_2, z_1, \ov z_2, t) \in Z^{-1}_{j,w,*}\times \mathcal I$.
Notice that the expression in the first logarithm never equals zero. Also, the expression in the second logarithm is well defined and never equal to zero since $(z_1,z_2)$ as well as
$(w_1,w_2)= (f_{j,\l_1} (z_1),f_{j,\l_2} (z_2))$ are pairings and thus their respective distance is at most $ \d_0\leq \frac14$.
Since, moreover, the set
$\La$ is simply connected, both logarithms in \eqref{40} are well defined and  we can and will take the principle branch since for $(\l_1 ,\ov\l_2, z_1, \ov z_2)=(\l ,\ov\l, z, \ov z) \in \De\cap Z^{-1}_{w,*}$ both expressions in the arguments of the logarithms are real positives. We thus have a properly defined  map $\Phi_j$ on every set $Z^{-1}_{j,w,*}$. 

The map $\Phi_j$ is in fact a  global well defined map on the union $\bigcup_{w , *} Z^{-1}_{j,w,*}\times \mathcal I$. In order to see this, consider  two sets $Z^{-1}_{j,w,*}$ and $Z^{-1}_{j,w',*'}$ having nonempty intersection. Then $\Delta \cap Z^{-1}_{j,w,*}\cap Z^{-1}_{j,w',*'}$ is a non-empty non-analytic subset of $Z^{-1}_{j,w,*}\cap Z^{-1}_{j,w',*'}$ and $\Phi_j$ restricted to $(\Delta \cap Z^{-1}_{j,w,*}\cap Z^{-1}_{j,w',*'})\times I$ is real and coincides with the given potential $\ph_j$. The map $\Phi_j$ is thus the desired extension of $\ph_j$
to $\mirror^{-1}\times\mathcal I$.

Given this extended potential and using the inclusion in \eqref{17}, we can now consider the extended operator $L_{j,\l_1,\ov \l_2, t}$ acting on functions $g\in \analytic$ by
\beq \label{19}
L_{j,\l_1,\ov \l_2, t}g (w_1,\ov w_2 ) = \sum _{z_1, z_2}\exp \left(\Phi_{j,\l_1,\ov \l _2 , t} (z_1, \ov z_2 )\right) g(z_1 , \ov z_2)
\eeq
where the summation is taken over all $1$--pairings $(z_1, z_2)$ such that $f_{j,\l_i}(z_i)=w_i$ , $i=1,2$. As for the initial real operator $\pf_j$ it is convenient to write simply $L_j$ instead of $L_{j,\l_1,\ov \l_2, t}$ when it is clear that the parameters $\l_1,\l_2,t$ are fixed.

  In the next proposition we will see that the image function $L_{\l_1,\ov \l_2, t}g\in \analytic$ provided the initial real operator $\pf_{\l , t }$ is bounded. This will allow us  to iterate the  operator and this will be done again in a non-autonomous way: (in \eqref{L^n} we use the abbreviated notation $L_k:=L_{k,\lambda_1,\ov\lambda_2,t}$)
\beq\label{L^n}
L_j^n = L_{j+n-1}\circ ...\circ L_j
\eeq
is the extension of $\pf_j^n$ defined in \eqref{non auto comp pf}.
Notice that the cocyle properties of inverse branches along with Lemma \ref{16} show that $L_j^n$ can also be defined by formula \eqref{19} if one replaces the potential $\Phi_{j,\l_1,\ov\l_2,t}$ by
$$
S_n \Phi_{j,\l_1,\ov\l_2,t}(z_1,\ov z_2)= \sum_{k=0}^{n-1} \Phi_{j+k, \l_1,\ov\l_2,t}(f^k_{\l_1}(z_1), \ov {f^k_{\l_2}(z_2)})\,
$$
and where the summation is taken over all $n$--pairings $(z_1,z_2)$
such that $f_{j,\l_i}^n (z_i)=w_i$, $i=1,2$.

\bprop\label{14} Suppose that the real operator $\pf_{j,\l , t }$ is uniformly bounded for $j\in \Z$, $\l\in \La$ and $t\in I$
and that $r=diam(\La ) /2$ is sufficiently small such that \eqref{41} holds for all $1$--pairings. 
Then there exist $a>0$ such that, with $$\mathcal I= \{x+iy\in \C \; ; \; x\in I \, , \; y\in ]-a,a[\},$$
the extended operator $L_{j,\l_1, \ov \l _2 , t}$ is a, uniformly for $j\in \Z$, $(\l_1, \ov \l_2, t)\in \La \times \ov \La \times \mathcal I$, bounded operator of $\analytic$. Moreover, if $\l_1=\l_2 =: \l$ and if $t\in I$ is real, then each operator $L_k:= L_{k,\l ,\ov \l , t}$ preserves $\ranalytic$ and there exists $K<\infty$ such that, for every function $h\in \analytic$,
\beq\label{26}
 \big| L^n_j h(w_1 , \ov w_2)-L^n_jh(w_1 , \ov w_1) \big |  \leq 
\pf_{j,\l ,t}^n \1 (w_1) \left( K+ \frac{ \ga_n^{-1}}{\d(1-\ga_n^{-1})} \right) \|h\| _\infty |w_1-w_2|
\eeq
where $(w_1, \ov w_2)\in \mirror _w$ and $n\geq 1$, and, as in \eqref{L^n} $L^n_j=L_{j+n-1}\circ\dots\circ L_j$.
\eprop
\noindent

\bpf  Let  $j\in \Z$, $(\l_1, \ov \l _2,w_1,\ov w_2)\in \mirror $ and let $t\in \mathcal I$ be complex.
For every $g\in \cA$ we have 
$$\left|L_{j,\l_1, \ov \l _2 , t} g (w_1,\ov w_2) \right|\leq \|g\|_\infty \sum_{z_1,z_2} \left| \exp \left( \Phi_{j,\l_1, \ov \l_2 , t}(z_1,\ov z_2)\right)\right|$$
where the summation is again over all corresponding $1$--pairings like in \eqref{19}. Therefore, it suffices to estimate the
series on the right hand side of this inequality in order to get a bound of the norm of the operator $L_{j,\l_1, \ov \l _2 , t}$ on $\analytic$.

 Now, if $(z_1, \ov z_2)$ be a $1$--pairing such that
$f_{j,\l_i}(z_i)=w_i$, $i=1,2$, then
\begin{align*}
\left|  \exp ( \Phi_{j,\l_1, \ov \l_2 ,t } (z_1, \ov z_2)) \right| &= 
 \left|  f'_{j,\l_1}(z_1) f'_{j,\l_2}(z_2)\right|^{- \frac{\Re t}{2}}
\exp\left\{  \frac{\Im t}{2} \arg \left( f'_{j,\l_1}(z_1)\ov {f'_{\l_2}(z_2)}\right)
\right\}\\
&\times \left|  \left(\frac{1+z_1\ov z_2}{1+w_1 \ov w_2}\right) ^{-t\frac{\tau}{2}} \right|\,.
\end{align*}
The choice of $r>0$ and \eqref{41} shows that $\left|  \arg \left( f'_{j,\l_1}(z_1)\ov {f'_{j,\l_2}(z_2)}\right)\right| \leq \tilde A$. Since $|\Im t |\leq a$ it follows that
$$
\exp\left\{  \frac{\Im t}{2} \arg \left( f'_{j,\l_1}(z_1)\ov {f'_{j,\l_2}(z_2)}\right)
\right\}\leq \exp \left\{ \frac a2 \tilde A\right\}.
$$
Clearly, $\left| \arg  \left(\frac{1+z_1\ov z_2}{1+w_1 \ov w_2}\right) \right| $ is bounded above uniformly with respect to $z_i, w_i$, $i=1,2$. Denote this bound again by $\tilde A$. Setting $B= \exp\left\{  a \tilde A  \frac{1+\tau }{2}
\right\}$ it follows that 
$$\left|  \exp ( \Phi_{j,\l_1, \ov \l_2 ,t } (z_1, \ov z_2)) \right| \leq B
 \left|  f'_{j,\l_1}(z_1) f'_{j,\l_2}(z_2)\right|^{- \frac{\Re t}{2}}
 \left|  \frac{1+z_1\ov z_2}{1+w_1 \ov w_2}\right| ^{-\frac{\tau}{2}\Re t}\,.
$$
An elementary calculation shows that there exists a constant $C<\infty$ independent of $z_i, w_i$, $i=1,2$, and $t\in I$, such that 
$$
 \left|  \frac{1+z_1\ov z_2}{1+w_1 \ov w_2}\right|^{-\frac{\tau}{2}\Re t} \leq C \sqrt{\frac{1+|z_1|^2}{1+|w_1|^2}\frac{1+|z_2|^2}{1+|w_2|^2}}^{-\frac{\tau}{2}\Re t}\,.
$$
Therefore, 
\begin{align*}
\left|  \exp ( \Phi_{j,\l_1, \ov \l_2 ,t } (z_1, \ov z_2)) \right|  \leq BC 
  | f'_{j,\l_1}(z_1)| ^{- \frac{\Re t}{2}}& \left(\frac{1+|z_1|^2}{1+|w_1|^2}\right)^{- \frac{\tau \Re t}{4}} \times \\
 &\times   | f'_{j,\l_2}(z_2)| ^{- \frac{\Re t}{2}} \left(\frac{1+|z_2|^2}{1+|w_2|^2}\right)^{- \frac{\tau \Re t}{4}} \, ,
\end{align*}
and thus the Cauchy-Schwarz inequality implies that
$$
 \sum_{z_1,z_2} \left| \exp \left( \Phi_{j,\l_1, \ov \l_2 , t}(z_1,\ov z_2)\right)\right|\leq BC 
\sqrt{\pf_{j,\l_1, \Re t}\1 (w_1)}\sqrt{ \pf_{j,\l_2, \Re t}\1 (w_2) }\,.
$$
By our assumptions there exists $M<\infty $ such that $\| \pf_{j,\l, t_0}\1\|_\infty \leq M$ for every $j\in \Z$, $\l \in \La $ and $t_0\in I$. This shows that 
\beq\label{61}
\| L_{j,\l_1, \ov \l _2 , t}\|_\infty \leq BC M\,.
\eeq
Suppose now that $\l_1=\l_2 =: \l$ and that $t\in I$ is real. In this case each operator $ L_{k,\l ,\ov \l , t}$ clearly preserves $\ranalytic$.

It remains to establish the distortion property \eqref{26}. We have
$$\big| L^n_j h(w_1 , \ov w_2)-L^n_jh(w_1 , \ov w_1) \big |  \leq I+II$$
where
\begin{align*}
I=\left| \sum\exp S_n\Phi_{j,\l ,\ov \l , t}(z_1, \ov z_1)\left(  h(z_1,\ov z_2) -h(z_1, \ov z_1)\right) \right|
\leq \pf_{j,\l ,t} ^n \1 (w_1)\, \nlip (h) \ga_n^{-1} |w_1-w_2|\,.
\end{align*}
Lemma \ref{39} gives an appropriate estimation for $\nlip (h)$ and thus
$$
I\leq \pf_{j,\l ,t} ^n \1(w_1) \, \frac{\|h\|_\infty}{\d(1-\ga_n^{-1})} \ga_n^{-1} |w_1-w_2|\,.
$$
The second term is equal to 
$$II=\left| \sum \left(\exp S_n\Phi_{j,\l ,\ov \l , t}(z_1, \ov z_2) -\exp S_n\Phi_{j,\l ,\ov \l , t}(z_1, \ov z_1)\right)h(z_1,\ov z_2)\right|  \,.$$
The following distortion estimate directly results from the complex version of Koebe's distortion theorem in the case $\tau =0$ and from Lemma 4.7 in \cite{MUmemoirs} if $\tau>0$:
$$\label{6}
\left| \frac{\exp S_n\Phi_{j,\l ,\ov \l , t}(f_{j,\l ,*}^{-n}(w_1),\ov{ f_{j,\l ,*}^{-n}(w_2)})}{\exp S_n\Phi_{j,\l ,\ov \l , t}(f_{j,\l ,*}^{-n}(w_1), \ov {f_{j,\l ,*}^{-n}(w_1)})} -1 \right| \leq K |w_1-w_2|\,\,  , \,\, w_1,w_2\in \D(w,\d) .
$$
Consequently,
$$II \leq \pf _{j,\l ,t} \1 (w_1 ) \| h\|_\infty K |w_1-w_2| $$
and, combining this estimate with the one of $I$, the desired Lipschitz property follows.
\epf

\section{Complexification of the invariant density}\label{complex}

We have to consider appropriate rescaled versions of the operators defined in the previous section. 
This section deals with the case where $\l_1=\l_2=:\l$ and $t\in I$ is real. Moreover, here and in the next section both parameters $\l , t$ are fixed and so we will frequently surpress them:
\beq\label{48}
\hat \pf_j = e^{-P_j (t)} \pf_j 
\quad\text{and} \quad
\hat L_j = e^{-P_j (t)} L_j 
\quad , \quad j\in \Z\,.
\eeq
The number $P_j (t)$ is usually called the topological pressure.
Assume that for these rescaled operators there exist strictly positive functions  $\hat \rho_j \in C_b^0 (U)$ such that, for some $M<\infty$ and for every $j\in \Z$ and $n\geq 1$,
\beq \label{24}
\| \hat \pf^n_j\|_{\infty } \leq M 
 \quad \text{and} \quad 
 \hat \pf^n_{-n+j} \1  \to \hat \rho _j
\,.
\eeq
where the limit is with respect to the sup-norm as $n\to \infty$.
Then clearly
$$\hat\pf_j \hat \rho_j =\hat \rho_{j+1}, \quad  j\in \Z\, ,$$
 and, for this reason, these functions are called invariant densities.
The aim now is to extend the invariant densities to holomorphic functions of $\cA_\R$ such that \eqref{24} still holds.

\bprop\label{25}
Suppose \eqref{24} does hold.  Then, for every $j\in \Z$, the sequence $\hat L^n_{-n+j}\1 $ converges uniformly on compact sets to some function of $\cA_\R$.
These limit functions are extensions of $\hat\rho_j$ and they will be denoted by the same symbol. 
Moreover, 
$$\left|\hat \rho (w_1, \ov w_2 ) - \hat \rho (w_1, \ov w_1) \right| \leq M(K+1) |w_1-w_2| \; , \;\; (w_1,w_2)\in \mirror_w\, ,$$
and the invariance property 
\beq \label{29}\hat L_j \hat \rho_j =\hat \rho_{j+1} \quad \text{holds on $ \mirror _w$ for every $j\in \Z$.} \eeq
\eprop
\smallskip

\bpf
Let $(w_1, \ov w_2)\in \mirror_w$.  The distortion property \eqref{26} implies that there exists $n_0\geq 0$ such that 
for every $n\geq n_0$
\beq\label{27}
\left| \hat L^{n}_{-n+j} \1 (w_1, \ov w_2 )- \hat L^{n}_{-n+j} \1 (w_1, \ov w_1 )\right| \leq \hat \pf ^n_{-n+j}\1 (w_1) (K+1) | w_1-w_2|
\,.
\eeq
Since  $ \hat L^{n}_{-n+j} \1 (w_1, \ov w_1 )=  \hat \pf ^n_{-n+j}\1 (w_1) $ it follows that
$$
\begin{array}{cl}
\left| \hat L^{n}_{-n+j} \1 (w_1, \ov w_2 )\right| &\leq \hat \pf ^n_{-n+j}\1 (w_1)
\Big( 1+ (K+1)|w_1-w_2|\Big)\\
&\leq M \Big( 1+ (K+1)\d \Big)\leq M(K+2)\quad \text{for every} \quad n\geq n_0\,.
\end{array}
$$
Therefore, the sequence $\left(\left| \hat L^{n}_{-n+j} \1 (w_1, \ov w_2 )\right|\right)_{n=0}^\infty$ is uniformly bounded above. Montel's Theorem thus applies and yields normality of the family $\left(\hat L^{n} _{-n+j}\1\right)_{n=0}^\infty$.
Since the limit of every converging subsequence coincides with $\hat \rho_j$ on the non-analytic set $\Delta_w$
the whole sequence $\left(\hat L^{n} _{-n+j}\1\right)_{n=0}^\infty$ converges to one and the same limit and this limit belongs to $\analytic_{\mathbb R}$.

The invariance property \eqref{29} holds since it holds on the non-analytic set $\Delta_w$.
Finally, the limit functions have the required Lipschitz property because of \eqref{24} and \eqref{27}.
\epf

\
The obvious modification of this proof, where $\1$ is replaced by an arbitrary element of $\analytic_{\mathbb R}$, also shows that the normalized extended operators and invariant densities are uniformly bounded above.
Whenever \eqref{24} holds we may assume, increasing $M$ if necessary, that
\beq\label{x.1}
\|\hat L_j^n \|_\infty\leq M \quad \text{and}\quad \|\hat\rho_j\|_\infty \leq M \;\; \text{for every } \; j\in \Z \; ,\;\; n\geq 0 \,.
\eeq
In the condition \eqref{x.1}, $\hat\rho_j$ is the extended density and the sup-norms are taken on the whole mirror $\mirror_w$.

\

In the sequel we will need a different normalization. Let $l:C^0_b (U)\to \R$ be a bounded functional.
It naturally acts on functions of $\analytic_\R$:
if $h\in\analytic_\R$ then $g(z)=h(z,\ov z)$, $z\in U$, defines a function $g\in C^0_b (U)$ and thus
we can define $$l(h):= l(g)\,.$$ In particular, $l(\hat \rho_j)$ is well defined 
regardless of whether $\hat \rho _j $ is understood as the initial function of $ C_b^0(U)$ 
or the extended function that belong to $\analytic_\R$.

The functional $l$ is assumed to be uniformly positive on the density functions meaning that there exists $a>0$ such that
\beq \label{33}
l(\hat \rho _j) \geq a \quad \text{ for every $j\in \Z$.}
\eeq
\bexam \label{x.3}
Fix any point $\xi \in U$ and consider the functional $l$ defined by $l(g)=g(\xi)$. Such a functional is uniformly positive on the functions $\hat \rho_j$ in the sense of \eqref{33} as soon as the system is mixing. This holds in particular for the transcendental random systems considered in \cite{MyUrb2014}. Lemma 5.5 of that paper shows that
there exists $n_0\geq 1$ and $a>0$ such that
\beq \label{32}
\npf ^n _j \1 (\xi )\geq a \quad \text{for every } \; n\geq n_0\,.
\eeq
\eexam

\medskip

\noindent
Consider then
\beq \label{31}
\rho_j = \frac{\hat\rho _j}{l(\hat\rho _j)}\quad , \quad j\in \Z\,.
\eeq
Clearly,
$$
\lim_{n\to \infty}\frac{\pf_{-n+j}^n\1}{l(\pf_{-n+j}^n\1)}= \lim_{n\to \infty}\frac{\hat \pf_{-n+j}^n\1}{l(\hat\pf_{-n+j}^n\1)}=\rho_j
$$
and, because of \eqref{29}, the extended invariant densities satisfy
\beq \label{35}
\frac{L^n_j (\rho_j )}{l(L^n_j(\rho_j))}= \rho_{j+n}\quad\text{for every $j\in \Z$ and $n\geq 1$.}
\eeq
It is henceforth natural to consider maps $\Psi_{n,j}$ defined by
\beq \label{56}
\Psi_{j}^n(g)=\frac{L^n_j (g )}{l(L^n_j(g))}=\frac{\hat L^n_j (g )}{l(\hat L^n_j(g))}\quad\text{for every $j\in \Z$ and $n\geq 1$.}
\eeq

\blem \label{57}
For every $j\in \Z$ and $n\geq 1$, the map $\Psi_{j}^n$ is well defined on the following neighborhood of $\rho_j$ in $\cA$:
$$
U_j:=\left\{g\in \cA : \;\; \| g-\rho_j\|_\infty < \frac{a}{2(\|l\|_\infty M)^2} \right\}\,.
$$
\elem

\bpf
For $g\in U_j$ we have to check that 
$$l(\hat L^n_j(g)) = l(\hat L^n_j(\rho_j))+l(\hat L^n_j(g-\rho_j))\neq 0\,.$$
Since
$$ l(\hat L^n_j(\rho_j))=\frac{ l(\hat L^n_j(\hat \rho_j))}{l(\hat \rho_j)}=\frac{ l(\hat \rho_{n+j})}{l(\hat \rho_j)}\, ,$$
since $l( \hat \rho _{n+j} ) \geq a$ by \eqref{33} and, since $l(\hat \rho_j )\leq \| l\|_\infty M$ by \eqref{x.1},
we have 
\beq  \label{58}
l\big( \hat L ^n (\rho_j )\big)\geq \frac{a}{\|l\|_\infty M}\,.
\eeq
On the other hand, if  $g\in U_j$ then $\|l(\hat L^n_j(g-\rho_j)) \|_\infty \leq \| l \|_\infty M \|g-\rho_j \| _\infty< \frac{a}{2\|l\|_\infty M}$.
Altogether we get $l(\hat L^n_j(g)) > \frac{a}{\|l\|_\infty M} -\frac{a}{2\|l\|_\infty M} = \frac{a}{2\|l\|_\infty M} >0$.
\epf

\smallskip

\section{Contraction}\label{contraction}
We shall exploit in detail the convergence of the normalized iterated operators under the assumption that there is a uniform speed of the convergence in \eqref{24}. Let us make this precise now (see also the condition  \eqref{uniform speed} in Theorem \ref{thm intro 2}).  We keep in this section the setting and notation of Section~\ref{complex} and assume again that  \eqref{24} and \eqref{33} hold. We also recall that \eqref{24} implies \eqref{x.1}.

We now fix $\d>0$  sufficiently small such that 
\beq\label{small delta}
\frac{M}{a}\|l\|_\infty\left(M \left(K+1\right)+Q\|l\|_\infty M\right)\d \leq \frac14\,,\eeq
where $Q=\frac{M}{a}(K+1)$.
Notice that diminishing $\d$ does not influence the involved constants since $M$ does not depend on $\d$ and the distortion constant  $K$ becomes even better if $\d$ is replaced by a smaller constant.

\

We shall formulate now the precise condition which we shall need in the sequel:

\

\noindent{\bf Uniform speed.} {\it There exist bounded linear functionals $\nu_j \in \analytic _\R ' $ and there exists a sequence $\om_n\to 0$ such that
\beq \label{28}
\| \hat \pf^n_{j} (h _{| \Delta_w} )-\nu_j (h) \hat \rho_{j+n} \|_{\infty, \Delta_w} \leq \om_n \|h\|_{\infty ,  \mirror_w} \quad \text{for every } \;\; h\in \analytic _\R \; , \;\; n\geq 1 \,.
\eeq
}

In order to avoid any confusion we indicated here the domain on which the sup-norm is taken. So on the left hand side of the inequality the supremum is taken over all points of the diagonal $\Delta_w$, which is identified with $U$, whereas 
on the right-hand side one takes into account the whole mirror $\mirror_w$.

We have chosen the notation $\nu_j$ since typical examples of these functionals are 
the measures of \eqref{conf measures} that often are called conformal measures.

\blem\label{x.2}
Assume that \eqref{24}, \eqref{33} and \eqref{28} hold. Then 
\beq \label{36}\nu_j (\hat \rho_j ) =1 \quad \text{, $\;\; j\in \Z$.}
\eeq
\elem

\bpf
Apply \eqref{28} with $h= \hat \rho _j$ and use the invariance property \eqref{29} in order to get
$$
\| \npf _j^n {\hat \rho_j } -\nu_j (\hat \rho_j )\hat \rho_{j+n}\|_{\infty, \Delta_w}
= | 1- \nu_j (\hat \rho_j )| \| \hat \rho_{j+n}\|_{\infty, \Delta_w}
\leq \om _n \|\hat \rho_j \|_{\infty ,  \mirror_w}
$$
By \eqref{x.1}, $ \|\hat \rho_j \|_{\infty ,  \mirror_w}\leq M$. On the other hand, 
\eqref{33} implies that $\| \hat \rho_{j+n}\|_{\infty, \Delta_w} \geq a/\| l \| _\infty$.
Since $\om_n\to 0$ as $n\to \infty$ we thus must have $\nu_j (\hat \rho_j ) =1$.
\epf

Let us now focus on $\pf ^n_0$, $n\geq1$, and use the simplified notations
$$
\pf^n = \pf ^n_0, \ L^n = L ^n_0,\ \nu =\nu_0,\  \rho = \rho_0, \ \hat \rho = \hat\rho_0,\ \Psi_n= \Psi_{n,0}. 
$$
Concerning the functional $l$, we already have explained the action of this functional on $\analytic _\R$.
It also can be extended to $\analytic$ by first extending it to complex functions in the usual way and then to functions $h\in \analytic$ by $l(h):= l( h_{|\Delta_w })$. 
Remember also the map  $\Psi_n$ given by
$ \Psi_n (g) = \frac{L^n (g)}{l(L^n (g))} $
is,  for every $n\geq 1$, well defined on the neighborhood $U_0$ of $\rho=\rho_0 $ (see  Lemma \ref{57}).

\bprop\label{34}
Suppose that \eqref{24}, \eqref{33} and the uniform speed condition hold.
Then, for every $\d\in ]0,\d_0]$ sufficiently small there exists $n\geq 1$ such that the differential of  $\Psi_n$ at $ \rho$ satisfies
$$\|   D_{\rho }\Psi^{n}\|_\infty\leq\frac{\sqrt{2}}{2}<1\, .$$
\eprop

\brem\label{34 rem}
The proof will show that the integer $n$ does not depend on the operators $\pf_j$ hence not on the functions $f_j$, $j\in Z$,
but only on the involved constants such as $a,M,\om_n$. In other words, $n$ is uniform for all families of operators as long as they satisfy the conditions \eqref{x.1}, \eqref{33} and the uniform speed with the same constants. This is in particular the case for 
all $\Psi_{n,j}$, $j\in \Z$.
\erem

\bpf
Let $h\in \analytic$.
From \eqref{35} we get $\Psi^n (\rho) = \rho_n$ and
\begin{align*}
\Psi^n (\rho +h) &=\frac{L^n(\rho)+ L^n (h)}{l(L^n(\rho ))+ l(L^n(h))}
=\frac{\rho_n+ L^n (h)/ l(L^n(\rho ))}{1+ l(L^n(h))/ l(L^n(\rho ))} \\
&= \rho_n+ \frac{L^n (h)}{ l(L^n(\rho ))} - \rho_n  \frac{l(L^n(h))}{ l(L^n(\rho ))} +o(\| h\| )\,.
\end{align*}
Hence,
$$D_{\rho }\Psi^n (h) = \frac{L^n (h)}{ l(L^n(\rho ))} - \rho_n  \frac{l(L^n(h))}{ l(L^n(\rho ))} \,.$$
Consider first the case where $h\in\analytic _\R$. It suffices to consider functions $h$ for which
$\|h\| _{\infty}\leq 1$.
If we evaluate the above expression at points $(w,\ov w)\in \Delta_w$ of the diagonal then we can use \eqref{28} and it follows that there are functions $\xi_n $ such that
$\|\xi_n \|_\infty \leq \om_n $ and such that
$$
\hat L^n (h) (w, \ov w) = \nu (h) \hat \rho _n (w) +\xi_n (w) \,.
$$
Consequently,
$$
\frac{L^n (h)}{ l(L^n(\rho ))}=\frac{\hat L^n (h)}{ l(\hat L^n(\rho ))}=\frac{\nu (h) \hat \rho _n +\xi_n }{l(\hat L^n(\rho ))}\quad \text{on} \quad \Delta_w\, .
$$
Thus
$$
D_{\rho }\Psi^n (h) _{|\Delta _w}=\frac{\nu (h) \hat \rho _n +\xi_n }{l(\hat L^n(\rho ))} 
-\rho_n \frac{\nu (h) l(\hat \rho _n )+l(\xi_n) }{l(\hat L^n(\rho ))}=
\frac{\xi_n -\rho_n l(\xi_n)}{l(\hat L^n(\rho ))} 
\,.
$$
This expression  can be estimated as follows.
From \eqref{58} we have
$l\big( \hat L ^n (\rho )\big)\geq \frac{a}{\|l\|_\infty M}$. For the same reasons, i.e. from \eqref{x.1} and \eqref{33}, we also have that
$ \|\rho _n\|_\infty=\frac{\| \hat \rho _n\| _\infty}{l(\hat \rho _n )}\leq \frac{M}{a} $. Altogether it follows that
\beq \label{37}
\| D_{\rho }\Psi^n (h) _{|\Delta _w}\|_\infty
\leq  \frac{\|\xi_n\|_\infty \left( 1+ \|\rho_n\|_\infty \|l\|_\infty\right)}{a/M\|l\|_\infty}
\leq \om_n \frac{M\|l\|_\infty}{a}\left( 1+\frac{M}{a}\| l \|_\infty\right)
\leq \frac14  \;.
\eeq
for all $n\geq n_0$ and some sufficiently large $n_0$.

For general points $(w_1,\ov w_2)\in \mirror$ we can proceed as follows. First of all we have
$$ D_{\rho }\Psi^n (h) (w_1, \ov w_2) = \frac{1}{l(\hat L^n (\rho ))}
\Big( \hat L^n (h)(w_1, \ov w_2) -\rho _n (w_1, \ov w_2) l(\hat L^n (h))\Big)\,.$$
We already have an appropriated estimate for the first factor. From the Lipschitz property  of $\hat \rho$ (Proposition \ref{25}) 
follows that
$$\left| \rho (w_1, \ov w_2) - \rho (w_1, \ov w_1) \right| \leq \frac{M(K+1)}{l(\hat \rho )}|w_1-w_2|
\leq  Q|w_1-w_2| \; , \quad (w_1, \ov w_2)\in \mirror_w\, $$
where, we remember, $Q=\frac{M}{a}(K+1)$.
If we combine this with the Lipschitz behavior of $L^nh$ given in \eqref{26} and use $|w_1-w_2|<\d$, we finally get
for large $n$
\begin{align*}
\Big|D_{\rho }\Psi^n (h) (w_1, \ov w_2) - & D_{\rho }\Psi^n (h) (w_1, \ov w_1)\Big| \leq
\frac{M}{a}\|l\|_\infty\left(M \left(K+\frac{8}{\d \ga _n}\right)+Q\|l\|_\infty M\right)\d\,.
\end{align*}
Remember now that $\d>0$ has been fixed small enough such that \eqref{small delta} holds.
This constant $\d$ being chosen, we can choose  $n$ sufficiently large such that $\frac{8}{\d \ga_n}\leq 1$.
Then
$$\Big|D_{\rho }\Psi^n (h) (w_1, \ov w_2) - D_{\rho }\Psi^n (h) (w_1, \ov w_1)\Big| \leq \frac14\,.$$
Combing this with \eqref{37}  implies that for real $h$ such that $\|h\|_\infty \leq 1$ we have, for this choice of $n$,
$$\|D_{\rho }\Psi^{n} (h) \|_\infty \leq \frac12\,.$$
 If $h\in\mathcal A$ is  arbitrary with $\|h\|_\infty=1$, then $h$ can be expressed as $h=h_1+ih_2$ where both $h_1$, $h_2$ are in $\mathcal A_{\R}$ and such that $ \max \{\|h_1\|_\infty, \|h_2\|_\infty \}\leq \|h\|_\infty =1$. 
  It suffices then to use the case of functions in $\mathcal A_{\mathbb R}$ 
of norm at most one in order to conclude this proof.
\epf

\section{Analyticity: the non-autonomous case}\label{analyticity non auto}

We now come to the final part where we investigate analytic dependence on the parameter $\l$. 
In this section we still continue with the non-autonomous case and thus with the notations introduced 
in the previous sections 3 to 6. The assumptions are also unchanged: \eqref{24}, thus \eqref{x.1}, \eqref{33}
and the uniform speed assumption \eqref{28} are kept throughout this section.

The first observation concerns the extended operators introduced in \eqref{19}.

\begin{prop}\label{59}
For every $j\in \Z$ and every $g\in \cA$, the map
$$ 
(t,\lambda_1, \ov{\lambda_2})\mapsto L_{j,\l_1,\ov \l_2, t}g \in \cA
$$
is holomorphic on $\mathcal I\times \Lambda\times \ov \Lambda$.
\end{prop}

\bpf

\fr Let $t_0\in \mathcal I$ and $\ep>0$ such that $\D(t_0,\ep)\subset \mathcal I$. We have to show that there are 
functions $h_{k_1,k_2,k_3}\in \cA$ such that for every $(t,\lambda_1, \ov{\lambda_2})\in\D(t_0,\ep)\times \Lambda\times \ov \Lambda$, we have a power series representation:
\beq\label{xw2}
L_{j,\l_1,\ov \l_2, t} g 
= \sum_{k_1,k_2,k_3\ge 0} h_{k_1,k_2,k_3} (\l_1-\l_0)^{k_1} (\ov \l_2-\ov \l_0)^{k_2} (t-t_0)^{k_3}.
\eeq
Every point $(w_1,\ov w_2)\in \mirror_w$ belongs to a disk $D_w$ for some $w\in U$. By Formula \eqref{xw1} we have well defined holomorphic functions 
$$
\La\times \ov \La\ni(\l_1,\ov\l_2) \mapsto 
{F_{j, \l_1, \ov \l_2, *}}_{\big |D_w\times \ov D_w}\in \cA_{\big |(D_w\times \ov D_w)\cap \mirror_w}
$$
ascribing to every $(w_1,\ov w_2) \in (D_w\times \ov D_w)\cap \mirror_w$ a $1$--pairing $(z_1,\ov z_2)$:
$$
F_{j, \l_1, \ov \l_2, *}(w_1,w_2) := (z_1,\ov z_2)
=(f_{j,\l_1,*}^{-1}(w_1), \ov {f_{j,\l_2,*}^{-1}(w_2)}).
$$
In consequence, the function
$$
\begin{aligned}
\D(t_0,\ep)\times\La\times \ov \La\ni&(\l_1,\ov\l_2) \longmapsto 
L_{j,\l_1,\ov \l_2, t} g_{\big |(D_w\times \ov D_w)\cap \mirror_w} =\\
&= \sum_*  \exp \Phi_{j,\l_1,\ov\l_2 , t} \circ F_{j, \l_1, \ov \l_2, *} g\circ F_{j, \l_1, \ov \l_2, *}
\in \cA_{\big |(D_w\times \ov D_w)\cap \mirror_w}
\end{aligned}
$$
is also holomorphic as the sum of an absolutely uniformly convergent series of holomorphic functions. Hence we have the representation,
$$
L_{j,\l_1,\ov \l_2, t} g|_{(D_w\times \ov D_w)\cap \mirror_w}
=\sum_{k_1,k_2,k_3\ge 0} h_{k_1,k_2,k_3;w} (\l_1-\l_0)^{k_1} (\ov \l_2-\ov \l_0)^{k_2} (t-t_0)^{k_3},
$$
where all the functions $h_{k_1,k_2,k_3;w}$ belong to $\cA|_{(D_w\times \ov D_w)\cap \mirror_w}$. From the uniqueness theorem for holomorphic functions, all these functions $h_{k_1,k_2,k_3;w}$, $w\in U$, glue to one element $h_{k_1,k_2,k_3}$ of $\cA$ giving rise to the representation \eqref{xw2}. The proof is complete.
\epf

Consider now a new Banach space $\mathcal A_\Z$ of all bounded sections $g=(g_j)_{j\in \Z}$
where $g_j\in \cA$ for every $j\in \Z$ and such that
$$\norm g \norm=\sup_{j\in \Z} ||g_j||_\infty.$$
The space $\mathcal A_\Z$ equipped with this norm $\norm \cdot \norm$ is a Banach space. 
One then considers the global operator $L_{\l_1,\ov\l_2,t}$ mapping $g=(g_j)_{j\in \Z}\in \analytic _\Z$
to the function $L_{\l_1,\ov\l_2,t}g\in \analytic _\Z$ which is defined by
$$
\left( L_{\l_1,\ov\l_2,t}g\right)_{j+1} = L_{j,\l_1,\ov\l_2,t}g_j \quad , \quad j\in \Z\,.
$$
In the same way, the map $\Psi^n_{j,\lambda_1,\ov{\l_2}, t}$ introduced in \eqref{56} gives rise to a global map $g\mapsto \Psi^n_{\lambda_1,\ov{\l}_2, t}(g)$ defined by
\beq\label{66}
(\Psi^n_{\lambda_1,\ov\l_2, t}(g))_{j+1} = \frac{L^n_{j,\l_1,\ov \l_2, t} (g_{j})}{l(L^n_{j,\l_1,\ov \l_2, t} (g_{j}))}
=\frac{\hat L^n_{j,\l_1,\ov \l_2, t} (g_{j})}{l(\hat L^n_{j,\l_1,\ov \l_2, t} (g_{j}))}
\quad , \quad j\in \Z .
\eeq
The integer $n\geq 1$ will be fixed such that the conclusion of Proposition \ref{34} holds.

Remember also that  for $t\in I$ real and for $\l=\l_1=\l_2$ the function 
$$\rho_{\l,\ov \l, t} = (\rho_{j,\l ,\ov \l, t})_{j\in \Z}$$ is a  fixed point of $\Psi^n_{\lambda, \ov{\l}, t}$ (see \eqref{35}).

\blem\label{60}
Let $\l_0\in \La$, let $t_0\in I$ be real and let $n\geq 1$. Then there exist $U_{\l_0 , t_0}$, an open neighborhood of 
$\rho_{\l_0,\ov \l_0, t_0}$ in $\mathcal A_\Z$ and an open neighborhood 
 $W_{\l_0 , t_0}$ of the point $(\l_0, \ov\l_0, t_0)$ in $\La \times \ov\La \times \mathcal I$   such that 
$ \Psi^n_{\lambda_1,\ov{\l_2}, t}$ is well defined on $U_{\l_0, t_0}$ for every $(\lambda_1,\ov{\l}_2, t)\in W_{\l_0 , t_0}$.
Moreover, the map
$$
U_{\l_0, t_0}\times W_{\l_0 , t_0}\ni (h,\lambda_1,\ov{\l}_2, t) \mapsto \Psi^n_{\lambda_1,\ov{\l}_2, t}(h)\in \mathcal A_\Z
$$
is holomorphic.
\elem

\bpf
First of all note that for every $j\in \Z$ and $n\geq 1$ the function  
$$
\mathcal A_\Z\times \La \times \ov\La \times \mathcal I\ni (h,\lambda_1,\ov{\l_2}, t) \mapsto L^n_{j,\l_1, \ov \l _2 , t}(h_j)\in \mathcal A
$$
is holomorphic since it is linear with respect to the first variable, holomorphic with respect to all three other variables
(Proposition \ref{59}), and one applies Hartogs' Theorem. Hence, also the function
\beq \label{xw3}
\mathcal A_\Z\times \La \times \ov\La \times \mathcal I\ni (h,\lambda_1,\ov{\l_2}, t) \mapsto l\(L^n_{j,\l_1, \ov \l _2 , t}(h_j)\)\in \C
\eeq
is holomorphic. Now, in order to conclude the proof, we shall 
find $U_{\l_0 , t_0}$, an open neighborhood of 
$\rho_{\l_0,\ov \l_0, t_0}$ in $\mathcal A_\Z $ and an open neighborhood 
 $W_{\l_0 , t_0}$ of the point $(\l_0, \ov\l_0, t_0)$ in $\La \times \ov\La \times \mathcal I$   such that 
$|l\(L^n_{j,\l_1, \ov \l _2 , t}(h_j)\)|$ is uniformly bounded below for every $h\in U_{\l_0, t_0}$ and for every $(\lambda_1,\ov{\l_2}, t)\in W_{\l_0 , t_0}$. This will tell us that all coordinates of the function $\Psi^n_{(\cdot, \cdot, \cdot)}(\cdot)$ are continuous and uniformly bounded, and ultimately  the function $\Psi^n_{\cdot, \cdot, \cdot}(\cdot)$ is holomorphic. 

\sp Let $n\geq 1$ be fixed. In order to find these neighborhoods we deduce from \eqref{x.1} that $ \| L^n_{j,\l_1, \ov \l _2 , t}1\|_\infty$ is uniformly bounded above with respect to $j\in\Z$ and $(\l_1, \ov \l _2 , t)\in \La \times \ov\La \times \mathcal I$.
Cauchy's Integral Formula thus implies that the map $(\l_1, \ov \l _2 , t)\mapsto  L^n_{j,\l_1, \ov \l _2 , t}\1$ is uniformly Lipschitz with respect to $j\in\Z$. Consequently, for every $\ep>0$ there exists a neighborhood $W_{\l_0 , t_0}$ of  $(\l_0, \ov\l_0,  t_0)$  such that for every $h\in \mathcal A_\Z$, we have that
\beq\label{62}
|L^n_{\l_1, \ov \l _2 , t}(h)- L^n_{\l_0, \ov \l _0 , t_0}(h)|=
\sup_{j\in \Z} \|  L^n_{j,\l_1, \ov \l _2 , t}(h_j)- L^n_{j,\l_0, \ov \l _0 , t_0}(h_j)\|_\infty \le\ep|h|\,.
\eeq
The existence of $U_{\l_0 , t_0}$  easily follows now from the above Lipschitz property  \eqref{62} along with  the estimate \eqref{58} of the proof of Lemma \ref{57}. 
\epf

We are now in a position to extend the invariant density $\rho_{\l_0,\ov \l_0 ,t_0}$ (i.e., to extend the function assigning the density  $\rho_{\l_0,\ov \l_0 ,t_0}$ to parameters $(\lambda_0,\ov\l_0,t_0)$) analytically to a neighbourhood of $(\l_0,\ov \l_0 ,t_0)$ by making use of the Implicit Function Theorem. Indeed, $\rho_{\l_0,\ov \l_0 ,t_0}$ is a fixed point of $\Psi^n_{\lambda_0,\ov{\l}_0, t_0}$,  the map $(h,\lambda_1,\ov{\l}_2, t) \mapsto \Psi^n_{\lambda_1,\ov{\l}_2, t}(h)$ is analytic (Lemma \ref{60}) and
Proposition~\ref{34} along with the Remark \ref{34 rem} imply that
$$
\norm D_{\rho_{\l_0,\ov \l_0 ,t_0}}\Psi^n_{\l_0,\ov \l_0 ,t_0}\norm 
= \sup_{j\in \Z}\| D_{\rho_{j,\l_0,\ov \l_0 ,t_0}}\Psi^n_{j,\l_0,\ov \l_0 ,t_0}\|_\infty
\leq \frac{\sqrt{2}}{2}<1
$$
provided $n$ has been chosen sufficiently large.
 In conclusion we get the following.

\bthm\label{64}
For every $(\l_0,t_0)\in \La\times I$ there exists  an open neighborhood $W_{\l_0 , t_0}$ in $ \La \times \ov\La \times \mathcal I$ of
 $(\l_0, \ov\l_0,  t_0)$, and $U_{\l_0 , t_0}$, an open neighborhood  of
$\rho_{\l_0,\ov \l_0, t_0} $ in $\cA_\Z$, 
along with an analytic map $(\l_1,\ov\l_2,t)\mapsto \rho_{\l_1,\ov\l_2,t}\in U_{\l_0 , t_0}$ such that 
$$
\Psi_{\l_1,\ov\l_2,t}(\rho_{\l_1,\ov\l_2,t}) =\rho_{\l_1,\ov\l_2,t} \quad \text{for every} \quad (\l_1,\ov\l_2,t)\in W_{\l_0 , t_0}\,.
$$
\ethm

\smallskip

\fr
Theorem \ref{thm intro 2} follows now easily.

\bpf[Proof of Theorem \ref{thm intro 2}]
An assumption of Theorem \ref{thm intro 2} is that there exists
$a>0$ and $z_0\in U$ such that $\hat \rho _{j,\l ,t} (z_0)\geq a$ for all $(j,\l ,t)$. 
This enables us to consider the functional $l:C^0_b (U)\to \R$ defined by $l(g):=g(z_0)$. It clearly satisfies \eqref{33}
and thus Theorem \ref{64} implies Theorem \ref{thm intro 2} provided the uniform speed condition \eqref{28} holds.
So, consider $h\in \analytic _\R$. By Lemma \ref{LipschitzLemma} the associated function $z\mapsto g(z)=h(z,\ov z)$
belongs to $Lip(U, \ka \d)$ with $\| g\|_{Lip, \ka \d}\leq C \|h\|_\infty$. It follows from the assumption  \eqref{uniform speed}
that there exists $\om_n\to 0$ such that
$$
\| \hat \pf^n_{j} g _{|\Delta_w} -\nu_j (g) \hat \rho_{j+n} \|_{\infty , \Delta_w} \leq \om_n \|g\|_{Lip, \ka\d}
\leq C\om_n \|h\|_{\infty  , \mirror_w}
$$
for every $n\geq 1$ and $ j\in \Z $.
This implies \eqref{28}  with $\om_n$ replaced by $C\om_n$.
\epf

\brem
Note that the uniqueness part of the Implicit  Function Theorem guarantees the functions $\rho_{\l,\ov\l,t}$, $t\in I$ being real, to coincide with the ones resulting from Proposition~\ref{25}.
\erem

\medskip

\section{Analyticity: the random case}

The final part of this paper is devoted to random dynamics. So we now consider the following setting. 
Let $X$ be an arbitrary set and $\mathcal B$ a $\sg$--algebra on $X$. We consider a complete probability space $(X,\mathcal B , m)$. As usual, the randomness will be modeled by an invertible map $\th : X\to X$ preserving the measure $m$.
All objects like functions and operators do now depend on $x\in X$ instead of the integer dependence $j\in \Z$ in the non-autonomous case. In particular, we consider functions $f_{x,\l}$, $x\in X$ and $\l\in \La$, that satisfy the conditions described in Section 2. In the random case one has to require in addition that these functions are measurable. This means that  the map 
$(x,z)\mapsto f_{x,\l}(z)$ is measurable for  every $\l\in \La$. We are interested in the dynamics of the random compositions
$$
f_{x,\l}^n = f_{\shift^{n-1} (x),\l}\circ ... \circ f_{x,\l}, \quad n\geq 1,$$
where $\l\in \La$ and $x\in X$. The associated radial Julia set $\J_r(f_{x,\l})$ is defined by the formula \eqref{51}
with functions $f_j, f_{j+1},...$ replaced by $f_{x,\l}, f_{\th (x),\l},...$.

The space of analytic functions $\cA_\Z$ is now replaced by $\cA_X$. It has the same meaning as before except that the functions depend measurably on $x\in X$. Thus, $g\in \cA_X$ if $(z_1,\ov z_2)\mapsto g_x(z_1,\ov z_2 )$ is holomorphic on $\mirror_w$ for every $x\in X$, if $x\mapsto g_x(z_1,\ov z_2 )$ is measurable for every $(z_1,\ov z_2 )\in \mirror_w$ and if 
$$
\norm g\norm := \ess\sup_{x\in X} \| g_x\|_\infty <\infty.
$$
The transfer operators $\pf_{x,\l,t}$ must also have measurable dependence on $x\in X$ in the sense that each function
$$
X\ni x\longmapsto \pf_{x,\l,t}g(z_1,\overline z_2)\in\C
$$
is measurable for all arguments $\l, t, g, (z_1,z_2)$ fixed in their appriopriate domains. Notice that one can show with the help of the Measurable Selection Theorem (see \cite{Cra02}) that this is indeed the case. In the case of transcendental functions this has been worked out in Lemma 3.6 of \cite{MyUrb2014}. In this case, the invariant densities $\rho_{x,\l ,t}$ as well as their extensions $\rho_{x,\l,\ov \l ,t}$ also depend measurably on $x\in X$ since they are obtained as a limit of measurable maps
(see \eqref{24} and Proposition \ref{25}). Clearly, exactly as for the above composition of the functions $f_{x,\l}$, the iterated operators are of the form $\pf_{x,\l,t}^n = \pf_{\th ^{n-1} (x),\l,t}\circ...\circ \pf_{x,\l,t}$. In the same way, the definitions given in the part on non-autonomous dynamics have straightforward counterparts. For example, 
the invariance of the density is the relation $\npf_{x,\l ,t} \hat \rho_{x,\l ,t}=\hat \rho_{\th (x) , \l ,t}$
and the uniform speed assumption \eqref{28} takes on the following form:
\beq \label{63}
\| \hat \pf^n_{x,\l,t} h -\nu_x (h) \hat \rho_{\shift ^n(x),\l,t} \|_{\infty , \Delta_w} \leq \om_n \|h\|_{\infty , \mirror_w} \quad \text{for every } \;\; h\in \analytic _\R \; , \;\; n\geq 1 \,.
\eeq
Also, the definition of the global map $g\mapsto \Psi_{\lambda_1,\ov{\l}_2, t}(g)$, $g\in \cA_X$, is 
$$
(\Psi_{\lambda_1,\ov\l_2, t}(g))_{\th (x)} = \frac{L_{x,\l_1,\ov \l_2, t} (g_{x})}{l(L_{x,\l_1,\ov \l_2, t} (g_{x}))}
\quad , \quad x\in X ,
$$
where again $l$ is a functional that satisfies \eqref{33}.
Proceeding now exactly as in the previous section and applying the Implicit Function Theorem in the Banach space $(\mathcal A_X , \norm . \norm )$ we see that Theorem \ref{64} holds also in the present random setting.

The results can now be summarized as follows. Assume again that the expanding property \eqref{20} is satisfied with $\gamma_n$ independent of $\lambda$, that this family is of bounded deformation (Definition \ref{15}) and the bounded distortion of the arguments of \eqref{41} holds. Finally, we assume that the, most natural in this context, thermodynamical formalism property \eqref{24} holds with some universal (i.e., independent of $\lambda$) constant $M$.

\bthm \label{thm analyticity}
Suppose the following:
\ben
\item There exists a bounded functional $l:C^0_b (U)\to \R$ that is uniformly positive on the invariant densities (see \eqref{33}).
\item The uniform speed condition \eqref{63} holds with some constants $\om_n$ independent of $\lambda$.
\een
Then, the map $(\l_1,\ov\l_2,t) \mapsto \rho_{\l_1,\ov\l_2,t}\in \mathcal A_X$ is analytic. In particular for a.e. $x\in X$ the map $(\l_1,\ov\l_2,t) \mapsto \rho_{x,\l_1,\ov\l_2,t}\in \mathcal A_X$ is analytic.
\ethm

\brem
Note that the uniqueness part of the Implicit Function Theorem guarantees the functions $\rho_{\l,\ov\l,t}$, $t\in I$ being real, to coincide with the ones resulting from Proposition~\ref{25}.
\erem

\brem \label{68}
In fact, in this theorem we also could include real analyticity of the expected pressure as defined in the transcendental case in \eqref{45} and established in Lemma \ref{46}.
\erem

\section{Transcendental random systems}\label{analyticity rtf} 
In this last part we apply the preceding results to the case of transcendental random systems.
Such systems have been considered in \cite{MyUrb2014} and the full thermodynamical formalism including spectral gap property has been established there. We here complete the picture in establishing analyticity in this  general context. As a consequence we get a proof for the particular example in the Introduction (Theorem \ref{thm intro 1}). 

Assume now that the functions $f_{x,\la}$ are transcendental functions and that this family consists of \emph{transcendental random systems} as defined in \cite{MyUrb2014}. 
We use notation from that paper such as $\cJ_{x,\l}$ for the Julia set of $(f_{x,\l}^n)_{n\geq 1}$. Clearly, the radial Julia set $\rad (f_{x,\l} )\subset \cJ_{x,\l}$. 
Here are some other notions from \cite{MyUrb2014} that are necessary for the present work.
First of all, the following mild technical conditions are used in \cite{MyUrb2014} with the same enumeration:

\setcounter{condition}{1}
\begin{condition}\label{C2}
There exists $T>0$ such that 
$$\Big(\jul _{x,\l}\cap \D_T \Big)\cap f_{x,\l}^{-1}\left(\jul_{\shift (x),\l} \cap \D_T\right) \neq \emptyset
\;\; , \;\; x\in X \text{ and } \l\in \La\,.$$
\end{condition}
\setcounter{condition}{3}
\begin{condition} \label{C4}
For every $R>0$ and $N\geq 1$ there exists $C_{R,N}$ such that
$$|\left( f_{x,\l}^N\right)'(z)|\leq C_{R,N} \quad \text{for all} \quad 
z\in \D_R\cap f_{x,\l}^{-N} \big( \D_R\big ) \; , \; x\in X \text{ and } \l\in \La\,.$$
\end{condition}

\noindent
Then, there must be some common bound for the growth of the (spherical) characteristic functions $\sT_{x,\l}(r)= \sT(f_{x,\l},r)$ of $f_{x,\l}$, $x\in X$ and $\l\in \La$.
We use here a stronger version of the Condition 1 in \cite{MyUrb2014} and would like to mention that this is only used in order to show that the expected pressure function has a zero (see Proposition \ref{t1ep2}):

\medskip

\noindent
{\bf Condition 1'.}
{\it
There exists $\rho>0$ and $ \iota >0$ such that 
\beq\label{49}
\iota r^{\rho } \leq \sT _{x,\l} (r) \leq \iota^{-1} r^{\rho } \quad \text{for all } \;
r\geq 1\;  , \; x\in X \text{ and all } \;\l\in \La\,.
\eeq
}

\bdfn\label{dfn hyperbolic}
The transcendental random family $(f_{x, \l})_{x\in X,\l\in\La}$ is called:
\ben
\item \emph{Topologically hyperbolic}
if there exists $0<\d_0\leq \frac14$ such that for every $x\in X$, $\l\in \La$, $n\geq 1$ and $w\in \jul_{\shift^n(x),\l}$
all holomorphic inverse branches of $f_{x,\l}^n$ are well defined on $\D (w, 2\d_0)$.
\item \emph{Expanding} if there exists $c>0$ and $\g >1$ such that 
$$|(f_{x,\l}^n)'(z)|\geq c\g ^n$$
for every $z\in \jul_{x,\l} \setminus f_{x,\l}^{-n} (\infty)$ and every $x\in X$, $\l\in \La$. 
\item \emph{Hyperbolic} if it is both topologically hyperbolic and expanding.
\een
\edfn

\bdfn \label{balanced growth}
The transcendental random family $(f_{x, \l})_{x\in X,\l\in\La}$ satisfies the balanced growth condition if
 there are $\al_2 > \max\{0, -\al_1\}$ and $\kappa \geq 1$ such that for every $(x,\l)\in X\times \La$ and 
every $z\in f_{x,\l}^{-1}(U )$,
 \begin{equation}\label{eq intro growth}
 \ka^{-1} \leq \frac{|f_{x,\l}'(z)|}{(1+|z|^2)^{\frac{\al _1}{2}}(1+|f_{x,\l}(z)|^2)^{\frac{\al_2}{2}}} \leq  \ka
\, . \end{equation}
\edfn

\noindent
In the following we always assume that the above conditions are satisfied. 

\bdfn\label{50}
A transcendental holomorphic random family $(f_{x, \l})_{x\in X,\l\in\La}$ will be called \emph{\whatever} if 
\begin{enumerate}

\sp\item the base map $\th:X\to X$ is ergodic with respect to the measure $m$,

\sp\item the system $(f_{x, \l})$ is hyperbolic, 

\sp\item the balanced growth condition is satisfied,

\sp\item the Conditions 1', 2 and 4 hold.
\end{enumerate}
\edfn

\noindent
In this context, the right potential to work with is $\ph_{\l , t}$ as defined in \eqref{47} but with $\tau = \al_1 +\tau'$
where $\tau' <\al_2$ is arbitrarily close to $\al_2$ such that 
\beq\label{x.11}
t>\rho/\tau > \rho /\a\quad , \quad  \a=\a_1 +\a_2\,.
\eeq
With such a choice, 
the following has been shown in \cite{MyUrb2014}:

\medskip

{\it
- The full thermodynamical formalism holds. In particular, there exist $\nu_{x,t}$, the Gibbs states, in fact generalized eigenmeasures of dual transfer operators, and unique equilibrium states 
$$
\mu_{x,t}=\hat\rho_{x,t} \nu_{x,t},\ \  \nu_{x,t}(\hat \rho_{x,t} )=1.
$$
Moreover, for every $t>\rho/\a$, there are constants $A_t,C_t<\infty$ and $\e_t>0$ such that 
\beq \label{54}
\hat\rho_{x,t} (z) \leq C_t (1+|z|)^{-\e_t t} \quad \text{and}\quad \| \hat\rho_{x,t}\|_\infty \leq A_t \quad 
\text{for all }  z\in U \; \text{and } \; x\in X\,.
\eeq
- The normalized iterated transfer operator converge exponentially fast (Theorem 5.1 (2)). 
}

\medskip

For \whatever  transcendental random families one has the bounded deformation property. Indeed, 
the following uniform control is a complete analogue of Lemma 9.7 in \cite{MUmemoirs} and can be shown with exactly the same normal family argument as in the proof given in  \cite{MUmemoirs}. 
Let us recall that $\La =\D(\l_0 , r)$.

\blem\label{lemma 9.7}
For every $\e >0$ there exists $0<r_\ep <r$ such that
$$ 
\left| \frac{f_\l '(f^{-1}_{\l , *}(w))}{f_{\l_0} '(f^{-1}_{\l_0 , *}(w))} -1\right|<\ep
$$
for every inverse branch $f^{-1}_{\l , *}$ defined on $\D(w,\d_0)$, $w\in U$, and every $\l \in \D(\l_0 , r_\ep )$.
\elem

\fr If we combine this with Koebe's Distortion Theorem (see for ex. Theorem 2.7 in \cite{McMullen94})
 then it follows that the condition \eqref{bdd def 2} of the bounded deformation property always holds. The first property of the bounded deformation property \eqref{bdd def 1} holds for many families
(see again \cite{MUmemoirs}) and clearly for the exponential family in Theorem \ref{thm intro 1}.

\subsection{Expected pressure}
Fix $t>\rho / \al$ and let us first discuss the  numbers $P_{x,\l}(t)$ of \eqref{48}. They depend on the transfer operator
$\pf_{x,\l}$ which itself has been defined with the auxiliary parameter $\tau$ such that \eqref{x.11} holds. Let us indicate for a moment this dependence by a superscript $\tau$: $\pf_{x,\l}^\tau$, $P_{x,\l}^\tau(t)$ and let $\nu_{x,\l}^\tau$ denote the associated Gibbs states (conformal measures) such that \eqref{conf measures} holds: $\pf_{x,\l}^{\tau \,*} \nu_{\th (x),\l}^\tau= e^{P_{x,\l}^\tau(t)}\nu_{x,\l}^\tau$.

If $\tau'\neq \tau$, for example if $\tau'=\tau+\Delta\tau >\tau$, then the potentials of the operators corresponding to $\tau, \tau'$ respectively are related by
$$
|f'_{x,\l}(z)|_{\tau '}^{-t} = |f'_{x,\l}(z)|_{\tau }^{-t}\frac{v(z)}{v(f_{x,\l }(z))}
\quad\text{where} \quad v(z)= (1+|z|^2)^{-\frac{t\Delta\tau}{2}}\,.
$$
Notice that this function $v$ does not depend on the parameter $x\in X$ and thus 
$$
\pf_{x,\l}^{\tau '} g = \frac {1}{v} \pf_{x,\l}^\tau \( v g\) \quad \text{,}\quad g\in C_b^0(U)\;.
$$
Then, $v\nu_{x,\l}^\tau$ is a finite measure and, with $\ga^{-1} _x=  v\nu_{x,\l}^\tau (\1 )= \int v d\nu_{x,\l}^\tau $, 
$\upsilon_{x,\l} =\ga _x  v\nu_{x,\l}^\tau$ a probability measure such that, for $g\in C_b^0(U)$,
$$
\begin{aligned}
\pf_{x,\l}^{\tau' \,*} \upsilon_{\th (x),\l}(g)&= \g_{\th (x)}\int \pf_{x,\l}^{\tau' }(g) v d\nu_{\th (x) , \l}^\tau
=\g_{\th (x)}\int \pf_{x,\l}^{\tau}(vg)  d\nu_{\th (x) , \l}^\tau\\
&= \g_{\th (x)}e^{P_{x,\l}^\tau(t)} \int vg  d\nu_{x,\l}^\tau = \frac{ \g_{\th (x)}}{\g_x}e^{P_{x,\l}^\tau(t)} 
\upsilon_{x,\l} (g)\,.
\end{aligned}
$$
Thus, for $\tau'$ we have Gibbs states $\nu_{x,\l}^{\tau '}=\upsilon_{x,\l}$ with corresponding pressures
\beq\label{xx.1}
P_{x,\l}^{\tau '}(t)= P_{x,\l}^\tau(t) + \log \g _{\th (x)} - \log \g _x \quad , \quad x\in X\,.
\eeq

\smallskip

Theorem~3.1 in \cite{MyUrb2014} states that $\sup_{x\in X} |P_{x,\l} (t) | <\infty$ for every $\l \in \La$.
This allow us now to introduce the expected pressure:
\beq \label{45}
\cE P _{\l } (t) =\int _X P_{x,\l} (t) dm(x) \,.
\eeq
The cohomological equation \eqref{xx.1} and invariance of $m$ implies that
$\cE P_\l (t)$ does not depend on the auxiliary parameter $\tau$.
The function $t\mapsto \cE P_\l (t)$ is well defined for $t> \rho/\al$.
Real analyticity of this function  is a consequence of the following result. Here and in the following $l$ is again a functional that satisfies 
 \eqref{33}. Notice that the existence of such a functional is guaranteed thanks to Example \ref{x.3}.

\blem\label{46}
For the expected pressure we have the following expression
$$\cE P_{\l } (t) = \int _X \log l\left( \pf _{x,\l ,t} \rho _{x,\l,t} \right) dm (x)$$
and the function $(\l,t)\mapsto \cE P_{\l } (t) $ is real analytic in $\La\times (\rho/\al , \infty)$.
\elem

\bpf
On the one hand we know that $\pf_{x,\l ,t} \hat \rho_{x,\l ,t} = e^{P_{x,\l }(t)}\hat \rho_{\th (x),\l ,t}$ and on the other hand $\pf_{x,\l ,t}  \rho_{x,\l ,t} = l\left( \pf _{x,\l ,t} \rho _{x,\l,t} \right)\rho_{\th (x),\l ,t}$. Since $ \rho_{x,\l ,t}  = \frac{\hat  \rho_{x,\l ,t} }{l(\hat  \rho_{x,\l ,t} )}$ it follows that
$$
\log \left( l\left( \pf _{x,\l ,t} \rho _{x,\l,t} \right) \right)=  P_{x,\l }(t) + \log \left( l(\hat \rho_{\th (x),\l ,t})\right) - \log \left( l(\hat \rho_{x,\l ,t})\right)\,.
$$
 It suffices to integrate this expression with respect to $m$ and to use that the measure $m$ is $\th$--invariant.
 The statement on analyticity results from this expression and the fact (see \eqref{xw3} and Theorem~\ref{thm analyticity}) that the function $(\l_1,\ov\l_2,t)\mapsto l\( L _{x,\l_1,\ov\l_2 ,t} \rho _{x,\l_1,\ov\l_2,t}\)\in\C$ is holomorphic.
\epf

\subsection{Bowen's Formula} This formula concerns a fixed random system or, in other words, a fixed parameter $\l\in \La$. We can therefore neglect this parameter throughout this subsection and consider a fixed random system $(f_{x})_{x\in X}$.
As our preparation for the proof of Bowen's Formula we are to deal with  expected pressure in greater detail.

\blem\label{52}
Let $t>\rho/\a$. Then for $m$-a.e. $x\in X$ and every $w\in \cJ_x$, 
$$
\EP(t)=\lim_{n\to\infty}\frac1n \log \pf_{\th^{-n}(x),t}^n\1(w).
$$
\elem

\begin{proof}
Taking $g_x:=\1$, item (2) of Theorem~5.1 in \cite{MyUrb2014} yields for every $n\ge 1$ that
$$
\big|\hat \pf_{\th^{-n}(x),t}^n\1(w)-\hat\rho_{x,t} (w) \big|\le B\vartheta^n
$$
for some $B\in(0,+\infty)$ and some $\vartheta\in (0,1)$. Since $\rho_{x,t} (w) >0$ this yields
$$
\left|\log\left( \frac{1}{\hat\rho_{x,t} (w) }\hat \pf_{\theta^{-n}(x),t}^n\1(w)\right)\right|\le \frac{B'}{\hat\rho_{x,t} (w) }\vartheta^n
$$
for every $n\ge 1$ with some constant $B'>0$. Using the standard Birkhoff's sum notation  $  S_n P_{y} = P_ y + P_{\th (y)}+...+P_{\th^{n-1} (y)}$, we have
$$\hat \pf_{\theta^{-n}(x),t}^n\1(w) = e^{-S_nP_{\th^{-n}(x)}(t)} \pf_{\theta^{-n}(x),t}^n\1(w) \,. $$
With this notation, it follows that
$$
\left|  \frac1n \log \pf_{\th^{-n}(x),t}^n\1(w) -\frac1n S_n P_{\th^{-n} (x)}(t) \right| \leq \frac{B'}{\hat\rho_{x,t} (w) }\frac{\vartheta^n}{n} + \frac{|\log (\hat\rho_{x,t} (w)) |}{n} \longrightarrow 0
$$
as $n\to\infty$.
The lemma now follows by applying Birkhoff's Ergodic Theorem to the function $x\mapsto P_x(t)$.
\end{proof}

\fr This characterization of expected pressure along with hyperbolicity of the system $(f_x)_{x\in X}$ and of Condition 1' allow us to establish the desired description of the behavior of the expected pressure.

\bprop\label{t1ep2}
The function $ t\mapsto\EP(t)$ is real-analytic (hence continuous) on $(\rho /\al , \infty)$,  strictly decreasing with 
$\frac{d}{d t} \EP (t) \leq -\log \ga <0$ and
satisfies
$$\lim_{t\downto \rho /\al }\EP(t)=+\infty \quad \text{and} \quad \lim_{t\to+\infty}\EP(t)=-\infty \,.$$
\eprop

\bpf
Analyticity has been established in Lemma \ref{46},
while the strict monotonicity and the limit at $+\infty$ are straightforward and standard whith the use of Lemma \ref{52}. The estimate of the derivative is due to the expanding property and the formula in Lemma \ref{52}. Here are the details:

Condition \ref{C2} implies that there exists $w_x\in \J _{x,\l} \cap \D_T$. Using the expanding property in Definition \ref{dfn hyperbolic} one can estimate as follows:
$$
\begin{aligned}
 \pf_{\th^{-n}(x),t+s}^n\1(w_x) &= \sum_{f^n_{\th^{-n}(x)}(z)=w_x} |f^n_{\th^{-n}(x)}(z)|_\tau ^{-t} \; | f^n_{\th^{-n}(x)}(z)|_\tau^{-s}\\
 &\leq (c\ga ^n)^{-s} (1+T^2)^{\frac{s\tau}{2}}  \pf_{\th^{-n}(x),t}^n\1(w_x) \text{ , } s>0\,.
 \end{aligned}
$$
Taking logarithms and dividing by $n$ yields
$$
\frac1n \log  \pf_{\th^{-n}(x),t+s}^n\1(w_x) -\frac1n \log  \pf_{\th^{-n}(x),t}^n\1(w_x)\leq 
\frac sn \log \left(\frac{(1+T^2)^{\frac{\tau}{2}}}{c}\right) -s \log \ga \,.
$$
The estimate of the derivative $\frac{d}{d t} \EP (t)$  follows now from differentiability of the expected pressure along with the formula in Lemma \ref{52}.

It remains to analyze the behavior of $\EP $ near $\rho /\al$. In order to do so, we will use Condition 1' along with Nevanlinna Theory as explained in \cite{MyUrb2014}. In the following we use the notations from that paper especially from the proof of Lemma 3.17. It is shown there that there exists $k>0$ and $\tilde R_0>0$ sufficiently large such that for every $R>\tilde R_0$ and every $w\in U\cap \D_R$
$$
\pf_x \1_{\D_R}  (w) \geq k R^{-(\al_2 -\tau)t}\int _{r_R}^R \frac{\sT _x (r)}{r^{\hat \tau t +1}}dr
$$ 
where $r_R= \om^{-1}(8\log R)$ and where $\om $ comes from Condition~1 in \cite{MyUrb2014}. This condition being replaced here by
Condition~1', we have $\om (r) = \iota r^\rho$ and $\sT _x (r)\geq \iota r^\rho$. Therefore, still with $\hat \tau = \al _1 + \tau$ and with $\hat k = k \iota$, we get, uniformly in $w\in U\cap \D_R$ and $ x\in X$, the lower bound
\begin{align*}
\pf_x\1_{\D_R}  (w) \geq & \hat k R^{-(\al_2 -\tau)t}\int _{r_R}^R \frac{dr}{r^{\hat \tau t -\rho+1}}\\
 =  &\hat k R^{-(\al_2 -\tau)t} \big(  \log R  - \log r_R  + O(\hat \tau t -\rho)\big) \,.
\end{align*} 
The number $\tau \in (0,\al_2)$ is chosen in dependence of $t$ arbitrarily close to $\al_2$ such that $t>\rho / (\al_1 +\tau ) > \rho / \al $ (see Remark~1.2 in  \cite{MyUrb2014}). It is therefore clear that for every $H>0$ one can choose $R=R_H>\tilde R_0$ and then
$t_H>\rho /\al$ such that for every $t\in (\rho/\al , t_H)$
$$\pf_x\1_{\D_R}  (w) \geq  H \quad \text{for every} \quad w\in U\cap \D_R \; , \; x\in X\,.$$
Now, if $\pf_x^{n-1} \1_{\D_R}  \geq H^{n-1}$ on $U\cap \D_R$ for some $n\geq 1$ then 
$$\pf_x^{n} \1 \geq \pf_x \left( \1_{\D_R} \pf_{\th (x)}^{n-1} \left(  \1_{\D_R} \right)  \right)\geq H^{n-1} \pf_x\1_{\D_R} \geq H^n \quad \text{on} \quad U\cap \D_R\,.$$
The formula $\lim_{t\downto \rho /\al }\EP(t)=+\infty$ follows now by induction and Lemma \ref{52}.
\epf

Now, let $\mu_{x,t}$ be the invariant family of measures defined in Section~\ref{intro}, i.e.,  $d\mu_{x,t}=\hat\rho_{x,t}d\nu_{x,t}$.

\blem\label{53} For every $t>\rho/\al$, 
 the function $(x,z)\mapsto\log|f_x'(z)|$ is $\mu_{x,t}$--integrable meaning that the integral 
$$
\chi_t:=\int_X\int_{J_x}\log|f_x'(z)|\,d\mu_{x,t}(z)\,dm(x)
$$
is well-defined and finite. Moreover, $\chi_t>0$. 
\elem

\brem \label{rem random measure}
The measures $(\mu_{x,t})_{x\in X}$ depend measurably on $x\in X$ and they are in fact disintegrations of a measure $\mu_t$ on the global space $\J = \bigcup_{x\in X}\{x\}\times \J_x$ having marginal $m$. Such a measure is often called random measure. Crauel's book \cite{Cra02} contains the general background related to random measures and \cite{MyUrb2014} all the details concerning the present setting. Also, Theorem 5.1 in \cite{MyUrb2014} tells us that 
$\mu_t$ is ergodic and invariant under the global skew product $(x,z)\mapsto (\th (x), f_x(z))$.
\erem

\bpf [Proof of Lemma~\ref{53}.] Let $t>\rho/\al$. The expanding property implies $\chi_t>0$. It remains to show that $\chi_t<\infty$.
It follows from the estimate given in \eqref{54} that 
$$\int_X\int _{\J_x} \log |z|\, d\mu_{x,t}\,dm(x) =\int_X\int _{\J_x} \log |z|\, \hat\rho_{x,t}\, d\nu_{x,t}\,dm(x) <\infty \; , \;\;\; x\in X\,,$$
and from invariance that
$$\int_X\int _{\J_x} \log |f_x (z) |\, d\mu_{x,t} \,dm(x)=\int_X\int _{\J_{\th (x)}} \log |z|\, d\mu_{\th (x),t}\,dm(x) <\infty \;, \;\;\; x\in X\,.$$
Thus, both functions $(x,z)\mapsto \log (1+|z|^2)$ and  $(x,z)\mapsto \log (1+|f_x (z)|^2)$ are $\mu_t$--integrable. From the balanced growth condition follows now $\mu_t$--integrability of the function $(x,z)\mapsto\log |f_x' (z)|$.
\epf

\fr Proposition \ref{t1ep2} yields the existence of a unique zero $h>\rho/\al $ of the expected pressure function. It turns out that this number coincides almost everywhere with the Hausdorff dimension of the radial Julia set.

\bthm[A version of Bowen's Formula]\label{Bowen Formula}\label{thm BF}
If $(f_x)_{x\in X}$ is an \whatever random  system, then
$$
\HD(\rad (f_x))=h \quad \text{for $m$-a.e. $x\in X$.}
$$
\ethm

\begin{proof}
Since $\mu_h$ is an ergodic measure, there is $M\in(0,+\infty)$ such that $$\mu_{x,h}\(J_r(x,M)\)=1 \quad \text{for all} \quad x\in X_1\, ,$$
 where $X_1\sbt X$ is some measurable set with $m(X_1)=1$, and
$$
J_r(x,M):=\big\{z\in J_r(x):\varliminf_{n\to\infty}|(f_x^n(z)|<M\big\}.
$$
First we shall prove that 
\beq\label{3abf1}
\HD(J_r(x,M))\ge h
\eeq
or $m$-a.e. $x\in X_1$. Fix $x\in X_1$ and $z\in J_r(x,M)$. Set
$
y:=(x,z)
$
and denote by $f_y^{-n}$ the inverse branch of $f_x^n$ defined on $\D (f_x^n(z),\d)$ mapping $f_x^n(z)$ back to $z$.
For every $r\in(0,\d)$ let $k:=k(y,r)$ be the largest integer $n\ge 0$ such that
\beq\label{2abf1}
\D (z,r)\sbt f_y^{-n}\(\D (f_x^n(z),\d)\).
\eeq
Since our system is expanding this inclusion holds for all $0\le n\le k$ and 
$$
\lim_{r\to 0}k(y,r)=+\infty.
$$
Fix $n=n_k\ge 0$ to be the largest integer in $\{0,1,2\ld,k\}$ such that 
$f_x^n(z)\in \D(0,M)$ and $s=s_k$ to be the least integer $\ge k+1$ such that $f_x^s(z)\in \D(0,M)$. 
It follows from Birkhoff's Ergodic Theorem that
\beq\label{1abf2.1}
\lim_{k\to\infty}\frac{s_k}{n_k}=1  
\eeq
for $m$-a.e. $x\in X_1$ , say $x\in X_2\sbt X_1$ with $m(X_2)=1$ and $\mu_{x,h}$-a.e. $z\in J_r(x,M)$, say $z\in J_r^1(x,M)$, with $\mu_{x,h}\( J_r^1(x,M)\)=1$. 
Since the random measure $\nu_h$ is $h$-conformal, 
i.e., since $\nu_{\th(x),h}(f_x(A))=\exp(P_x(h))\int_A|f'_x|^h_\tau d\nu_{x,h}$ for every Borel set $A$ on which $f_x$ is injective,
we get from \eqref{2abf1} and the definition of $n$ that
\beq\label{1abf1}
\nu_{x,h}(\D(z,r))
\le \nu_{x,h}\(f_y^{-n}\(\D(f_x^n(z),\d)\)\)
\le K_{z,M}^h\big|\(f_x^n)'(z)\big|^{-h} e^{-S_nP_x(h)}
,
\eeq
where the constant $K_{z,M}$ compensates the replacement of the $\tau$-derivative $|\(f_x^n)'(z)\big|_\tau$ by the Euclidean derivative $|\(f_x^n)'(z)\big|$. On the other hand
$\D(z,r)\not\sbt f_y^{-s}\(\D(f_x^s(z),\d)\)$. But since, by $\frac14$-Koebe's Distortion Theorem,
$$
f_y^{-s}\(\D(f_x^s(z),\d)\)\spt \D\(z,\frac14|(f_x^s)'(z)|^{-1}\d\),
$$
we thus get that $r\ge \frac14|(f_x^s)'(z)|^{-1}\d$. Equivalently,
$$
|(f_x^s)'(z)|^{-1}\le 4\d^{-1}r.
$$
By inserting this into \eqref{1abf1} and using also the Chain Rule, we obtain 
$$
\nu_{x,h}(\D(z,r))
\le (4K_{z,M}\d^{-1})^hr^he^{-S_nP_x(h)}\big|\(f_{\th^n(x)}^{s-n}\)'\(f_x^n(z)\)\big|^h.
$$
Equivalently:
\beq\label{1abf2}
\frac{\log\nu_{x,h}(\D(z,r))}{\log r}
\ge h+\frac{h\log(4K_{z,M}\d^{-1})}{\log r}-\frac{S_nP_x(h)}{\log r}+h\frac{\log\big|\(f_{\th^n(x)}^{s-n}\)'\(f_x^n(z)\)\big|}
{\log r}.
\eeq
Now, Koebe's Distortion Theorem yields
$$
f_y^{-n}\(\D(f_x^n(z),\d)\) \sbt \D\(z,K\d|(f_x^n)'(z)|^{-1}\).
$$
Along with \eqref{2abf1} this yields $r\le K\d|(f_x^n)'(z)|^{-1}$. Equivalently:
\beq\label{2abf2}
-\log r\ge -\log(K\d)+\log|(f_x^n)'(z)|.
\eeq
By Lemma~\ref{53} the function $(x,z)\mapsto\log|f_x'(z)|$ is $\mu_h$--integrable with $
\chi_h>0$. Therefore,
there exists a measurable set $X_3\sbt X_2$ with $m(X_3)=1$ and for every $x\in X_3$ there exists a measurable set $J_r^2(x,M)\sbt J_r^1(x,M)$ such that $\mu_{x,h}\(J_r^2(x,M)\)=1$ and 
\beq\label{3abf2}
\lim_{j\to\infty}\frac1j\log|(f_x^j)'(z)|=\chi_h\in(0,+\infty)
\eeq
for every $x\in X_3$ and every $z\in J_r^2(x,M)$, the equality holding because of 
Birkhoff's Ergodic Theorem. This formula, along with \eqref{1abf2.1} also yields
\beq\label{1abf3}
\lim_{n\to\infty}\frac1n\log\big|\(f_{\th^n(x)}^{s-n}\)'\(f_x^n(z)\)\big| =0
\eeq
for every $x\in X_3$ and every $z\in J_r^2(x,M)$. Since $\int_X P_x(h)\,dm(x)=0$, 
Birkhoff's Ergodic Theorem gives:
\beq\label{3abf3}
\lim_{j\to\infty}\frac1j S_jP_x(h)=0,
\eeq
for all $x\in X_4\sbt X_3$, where $X_4$ is some measurable set with $m(X_4)=1$. By combining this formula taken together with the three formulas \eqref{1abf3}, \eqref{3abf2}, and \eqref{2abf2}, and formula \eqref{1abf2}, we get
$$
\varliminf_{r\to 0}\frac{\log\nu_{x,h}(\D(z,r))}{\log r}\ge h
$$
for every $x\in X_4$ and every $z\in J_r^2(x,M)$. Since $\mu_{x,h}\(J_r^2(x,M)\)=1$, we thus obtain, using a version of Frostman's lemma (see, e.g., \cite{PUbook}, Theorem 8.6.3):
\beq\label{2abf3}
\HD(J_r(x))
\ge \HD(\mu_{x,h})\ge h
\eeq
for every $x\in X_4$ (with $m(X_4)=1$).

\sp We now shall establish the opposite inequality. We know from Lemma 3.19 in \cite{MyUrb2014} that for any $n\ge 1$ large enough, say $n\ge q\ge 1$,
$$
Q_n:=\inf\big\{\nu_{x,h}(\D(w,\d)):x\in X,\, w\in J_x\cap \D(0,n)\big\}>0.
$$
By the very definition of $J_r(x)$ we have that
\beq\label{6abf4}
J_r(x)=\bu_{n=q}^\infty J_r(x,n).
\eeq
Fix $n\ge q$. Keep both $x\in X_4$ and $z\in J_r(x,n)$ fixed (still $y:=(x,z)$), and consider an arbitrary integer $l\ge 0$ such that
\beq\label{1abf4}
f_x^l(z)\in \D(0,n).
\eeq
Let $r_l>0$ be the least radius such that
\beq\label{2abf4}
f_y^{-l}\(\D(f_x^l(z),\d)\) \sbt \D(z,r_l).
\eeq
But, by Koebe's Distortion Theorem, $f_y^{-l}\(\D(f_x^l(z),\d)\) \sbt \D\(z,K\d|(f_x^l)'(z)|^{-1}\)$; hence
\beq\label{5abf4}
r_l\le K\d|(f_x^l)'(z)|^{-1}.
\eeq
Formula \eqref{2abf4} along with Koebe's Distortion Theorem and \eqref{5abf4}, yield
\beq\label{4abf4}
\begin{aligned}
\nu_{x,h}(\D(z,r_l))
&\ge \nu_{x,h}\(f_y^{-l}\(\D(f_x^l(z),\d)\) \\
&\ge K_{z,M}^{-h}\big|\(f_x^l)'(z)\big|^{-h}e^{-S_lP_x(h)}\nu_{h,\th^l(x)}
    \(\D(f_x^l(z),\d)\) \\
&\ge K_{z,M}^{-h}Q_ne^{-S_lP_x(h)}\big|\(f_x^l)'(z)\big|^{-h} \\
&\ge (K\d K_{z,M})^{-h}Q_ne^{-S_lP_x(h)}r_l^h.
\end{aligned}
\eeq
where the constant $K_{z,M}$ again compensates the replacement of the $\tau$-derivative $|\(f_x^l)'(z)\big|_\tau$ by the Euclidean derivative $|\(f_x^l)'(z)\big|$. Therefore,
\beq\label{1abf5}
\frac{\log\nu_{x,h}(\D(z,r_l))}{\log r_l}
\le h-\frac{h\log(K\d K_{z,M})}{\log r_l}-\frac{S_{l}P_x(h)}{\log r_l}-\frac{Q_n}{\log r_l}.
\eeq
Formula \eqref{5abf4} equivalently means that
\beq\label{2abf5}
-\log r_l
\ge \log\big|\(f_x^l)'(z)\big|-\log(K\d)
\ge \hat\chi l-\log(K\d)
\eeq
with some $\hat\chi>0$ resulting from uniform expanding property of the system $(f_x)_{x\in X}$. Since the set of all integers $l\ge 1$ for which \eqref{1abf4} holds is infinite (as $z\in J_r(x,n)$), taking the limit of the right-hand side of \eqref{1abf5} over all such $l$s. and applying \eqref{2abf5}, \eqref{3abf3}, and also recalling that, by Birkhoff's Ergodic Theorem,
$$
\lim_{j\to\infty}\frac1jS_jP_x(h)=0,
$$
we obtain
$$
\varliminf_{r\to 0}\frac{\log\nu_{x,h}(\D(z,r))}{\log r}
\le \varliminf_{l\to \infty}\frac{\log\nu_{x,h}(\D(z,r_l))}{\log r_l}
\le h
$$
Consequently, $\HD(J_r(x,n))\le h$ for all $x\in X_4$. Together with \eqref{6abf4} and $\sg$-stability of Hausdorff dimension, we thus get that $\HD(J_r(x))\le h$ for all $x\in X_4$. Along with \eqref{2abf3} this finishes the proof.
\end{proof}

\

\subsection{Conclusion} All in all we now get the following analyticity result for the dimension of the radial limit set.

\bthm\label{thm main}
Suppose that the transcendental holomorphic random family 
$(f_{x, \l})_{x,\l}$ is \whatever
and let $h_\l$ be the fiberwise Hausdorff dimension of the radial limit set of $(f_{x, \l})_{x\in X}$, $\la\in \La$.
Then, $\l \mapsto h_\l $ is real-analytic.
\ethm

\bpf Bowen's Formula shows that $h_\l$ is the unique zero of the expected pressure function. The later is analytic and $\frac{\partial}{\partial t} \EP_\l (t) <0$ (Proposition \ref{t1ep2}). Therefore
the Implicit Function Theorem applies and yields analyticity of $\l \mapsto h_\l $.
\epf

\

It remains to discuss the initial example given in the Introduction.

\medskip

\bpf[Proof of Theorem \ref{thm intro 1}]
 Let $U=\{z\in\C: \Re z >1\}$. It is well known that $f_\eta=\eta e^z$ is a hyperbolic exponential map if $\eta $ is real and $\frac1{6e} < \eta < \frac{5}{6e}$. Moreover, the closure  $\ov{f_\eta ^{-1} (U) } \subset U$. An elementary calculation shows that there exists $b>0$ such that $\ov{f_\eta ^{-1} (U) } \subset U$ for every $\eta \in \Omega_b$
where
$$\Omega_b =\lt\{ \eta \in \C\; ; \;\; \frac1{6e} < \Re (\eta ) < \frac{5}{6e} \text{ and } | \Im (\eta ) | < b \rt\}\,.$$
It follows that $f_{\eta_n}\circ ...\circ f_{\eta_1}$, $n\geq 1$, defines an expanding non-autonomous sequence that satisfies
\eqref{20} for any choice of $\eta_1, \eta _2 ,... \in \Omega_b$. It is straightforward to see that we thus have for these parameters 
a \whatever transcendental random family provided that we explain the random model.

 In order to do so, let
$X= \D(0,1)^\Z$, $\mathcal B$ the Borel $\sg$-algebra, $m$ the infinite product measure of the normalized Lebesgue measure of the unit disk and $\th$ the left-shift map on $X$.

Consider now parameters $(a,r)$ such that $\D(a,r)\subset \Omega_{b/2}$. Let $x\in X$ and $x_0$ the $0$--coordinate of $x$.
We associate to these parameters the function
$\eta e^z = (a+ r x_0) e^z$. In such a way we get for every $x\in X$ a family $(a,r)\mapsto f_{\eta}$. However, this family only depends real analytically on $(a,r)\in \R^2$. In order to turn this into a holomorphic family it suffices to replace these parameters by complex ones with small imaginary part such that $a+ r x_0\in \Omega_b$ for every $x_0\in \D(0,1)$.
Theorem \ref{thm main} applies to this family.
\epf


\bibliographystyle{plain}

\end{document}